\documentclass[12pt,draftcls,onecolumn]{IEEEtran}
\usepackage{amssymb,amsmath,bm}
\usepackage[dvips]{graphicx}

\newtheorem{dl}{Theorem}
\newtheorem{tl}{Corollary}
\newtheorem{yl}{Lemma}
\newtheorem{xz}{Property}
\newtheorem{dy}{Definition}
\newtheorem{lz}{Example}

\newcommand{\R}{\mathbb{R}}
\newcommand{\C}{\mathbb{C}}
\newcommand{\n}{\mathcal{I}_n}
\newcommand{\m}{\mathcal{I}_m}
\newcommand{\yi}{\mathbf{1}}
\newcommand{\N}{\mathcal{N}}

\DeclareMathOperator{\sign}{sign} \DeclareMathOperator{\spans}{span}
\DeclareMathOperator{\sig}{sig} 
\DeclareMathOperator{\diag}{diag}\DeclareMathOperator{\myarg}{arg}
\newcommand{\sgn}[2]{\sig\left({#1}\right)^{#2}}
\newcommand{\A}{\boldsymbol{A}}
\newcommand{\B}{\boldsymbol{B}}
\newcommand{\D}{\boldsymbol{D}}
\newcommand{\T}{\boldsymbol{T}}
\newcommand{\x}{\boldsymbol{x}}
\newcommand{\y}{\boldsymbol{y}}

\newcommand{\w}{\boldsymbol{\omega}}
\newcommand{\de}{\boldsymbol{\delta}}
\newcommand{\al}{\boldsymbol{\alpha}}

\newcommand{\bxi}{\boldsymbol{\xi}}
\newcommand{\bI}{\boldsymbol{I}}
\newcommand{\G}{\mathcal{G}}
\newcommand{\bN}{\mathcal{N}}
\newcommand{\bb}{\boldsymbol{b}}
\newcommand{\zt}{\boldsymbol{\zeta}}

\begin{document}
\title{Finite-Time Consensus Problems for Networks of Dynamic Agents}

\author{Long~Wang~and~Feng~Xiao
\thanks{This work was supported by NSFC (60674050 and 60528007), National 973 Program (2002CB312200), and 11-5 project (A2120061303).}
\thanks{The authors are with the Intelligent Control Laboratory, Center for Systems and Control,
Department of Industrial Engineering and Management, and Department
of Mechanics and Space Technologies, College of Engineering, Peking
University, Beijing 100871, China~(e-mail: longwang@pku.edu.cn;
fengxiao@pku.edu.cn).}}

\maketitle

\begin{abstract}
In this paper, finite-time state consensus problems for
continuous-time multi-agent systems are discussed, and two
distributive protocols, which ensure that the states of agents reach
an agreement in a finite time, are presented. By employing the
method of finite-time Lyapunov functions, we derive conditions that
guarantee the two protocols to solve the finite-time consensus
problems respectively. Moreover, one of the two protocols solves the
finite-time weighted-average consensus problem and can be
successively applied to the systems with switching topology. Upper
bounds of convergence times are also established. Simulations are
presented to show the effectiveness of our results.
\end{abstract}


\IEEEpeerreviewmaketitle
\section{Introduction}
The theory of consensus or agreement problems for multi-agent
systems has emerged as a challenging new area of research in recent
years. It is a basic and fundamental problem in decentralized
control of networks of dynamic agents and has attracted great
attention of researchers. This is partly due to its broad
applications in cooperative control of unmanned air vehicles,
formation control of mobile robots, control of communication
networks, design of sensor networks, flocking of social insects,
swarm-based computing, etc.

In \cite{T. Vicsek A. Czirok E. Ben Jacob I. Cohen and O. Schochet},
Vicsek et al. proposed a simple but interesting discrete-time model
of finite agents all moving in the plane. Each agent's motion is
updated using a local rule based on its own state and the states of
its neighbors. The Vicsek model can be viewed as a special case of a
computer model mimicking animal aggregation proposed in \cite{C.
Reynolds} for the computer animation industry. By using graph theory
and nonnegative matrix theory, Jadbabaie et al. provided a
theoretical explanation of the consensus property of the Vicsek
model in \cite{A. Jadbabaie J. Lin and A. S. Morse}, where each
agent's set of neighbors changes with time as system evolves. The
typical continuous-time model was proposed by Olfati-Saber and
Murray in \cite{R. Olfati-Saber and R. M. Murray 1}, where the
concepts of solvability of consensus problems and consensus
protocols were first introduced. The authors used a directed graph
to model the communication topology among agents and studied three
consensus problems, namely, directed networks with fixed topology,
directed networks with switching topology, and undirected networks
with communication time-delays and fixed topology. And it was
assumed that the directed topology is balanced and strongly
connected. In \cite{W. Ren and R. W. Beard}, Ren and Beard  extended
the results of \cite{A. Jadbabaie J. Lin and A. S. Morse,R.
Olfati-Saber and R. M.  Murray 1} and presented mathematically
weaker conditions for state consensus under dynamically changing
directed interaction topology. In the past few years, consensus
problems of multi-agent systems have been developing fast and many
research topics have been addressed, such as agreement over random
networks \cite{Yuko,A. V. Savkin}, asynchronous information
consensus \cite{L. Fang}, dynamic consensus \cite{D. P. Spanos
dynamic consensus}, networks with nonlinear consensus protocols
\cite{L. Moreau 1}, consensus filters \cite{Reza Olfati-Saber2}, and
networks with communication time-delays \cite{R. Olfati-Saber and R.
M. Murray 1,Tanner,feng xiao,Dongjun Lee}. For details, see the
survey \cite{W. Ren survey} and references therein.

By long-time observation of animal aggregations, such as schools of
fish, flocks of birds, groups of bees, and swarms of social
bacteria, it is believed that simple, local motion coordination
rules at the individual level can result in remarkable and complex
intelligent behavior at the group level. We call those local motion
coordination rules {\it distributive protocols}. In the study of
consensus problems, they are called {\it consensus protocols}.

In the analysis of consensus problems, convergence rate is an
important performance index of the proposed consensus protocol. In
\cite{R. Olfati-Saber and R. M. Murray 1}, a linear consensus
protocol was given and it was shown that  the second smallest
eigenvalue of interaction graph Laplacian, called algebraic
connectivity of graph, quantifies the convergence speed   of the
consensus algorithm. In \cite{Yoonsoo Kim}, Kim and Mesbahi
considered  the problem of finding the best vertex positional
configuration so that the second smallest eigenvalue of the
associated graph Laplacian is maximized, where the weight for an
edge between two vertices was assumed to be a function of the
distance between the two corresponding  agents. In \cite{Lin Xiao},
Xiao and Boyd considered and solved the problem of the weight design
by using semi-definite convex programming, so that algebraic
connectivity is increased. If the communication topology is a
small-world network, it was shown that large algebraic connectivity
can be obtained \cite{R. Olfati-Saber-smallworld}. Although by
maximizing the second smallest eigenvalue of  interaction graph
Laplacian, we can get better convergence rate of the linear protocol
proposed in \cite{R. Olfati-Saber and R. M.  Murray 1}, the state
consensus can never occur in a finite time. In some practical
situations, it is required that the consensus be reached in a finite
time. Therefore, finite-time consensus is more appealing and there
are a number of settings where finite-time convergence is desirable.
The main contribution of  this paper is to address finite-time
consensus problems and discuss two effective distributed protocols
that can solve consensus problems in finite times. Furthermore, the
method used in this paper is of interest itself, which is partly
motivated by the work of \cite{V. T. Haimo}, in which continuous
finite-time differential equations were introduced as fast accurate
controllers for dynamical systems, and partly by the results of
finite-time stability of homogeneous systems in \cite{S. J. Bhat
ACC}.

The proposed protocols in this paper are {\it continuous} state
feedbacks, but they do not satisfy the Lipschitz condition at the
agreement states, which is the least requirement for finite-time
consensus problems because Lipschitz continuity can only lead to
asymptotical convergence. To prove that the two protocols
effectively solve finite-time consensus problems, we adopt the
theory of finite-time Lyapunov stability \cite{V. T. Haimo,S. J.
Bhat ACC}. Although for systems under general communication
topology, it is difficult and even impossible to find the valid
Lyapunov functions, we have succeeded in reducing the general case
into several special cases, in which appropriate Lyapunov functions
can be found. This paper shows that the convergence time is closely
related to the underlying communication topology, especially, the
algebraic connectivity for the undirected case. Large algebraic
connectivity can greatly reduce the convergence time. We also
compare the convergence rates of two systems under the same protocol
but with different protocol parameters and draw the conclusion that
one converges faster when agents' states differ a lot while the
other converges faster when agents' states differs a little. This
conclusion encompasses the case when one of the two systems adopts
the linear consensus protocol presented in \cite{R. Olfati-Saber and
R. M. Murray 1}. Therefore, in order to get shorter convergence
time, we can change those adjustable parameters in the proposed
protocols according to agents' states. By the same method, we also
study the case where the topology is dynamically changing.

This paper is organized as follows. Section II presents the
preliminary results in algebraic graph theory.  Section III sets the
basic setting in which the problem is formulated. The main
theoretical results are established in Section IV  and the
simulation results are given in Section V. Finally, concluding
remarks are provided in Section VI.

\section{Preliminaries}
In this section, we list some basic definitions and results in
algebraic graph theory. More comprehensive discussions can be found
in \cite{C. Godsil and G. Royal}.

Directed graphs will be used to model the communication topologies
among agents. A  {\it directed graph} $\mathcal{G}$ of order $n$
consists of a vertex set $\mathcal{V(G)}=\{v_i:i=1,2,\dots,n\}$ and
an edge set
$\mathcal{E}(\G)\subset\{(v_i,v_j):v_i,v_j\in\mathcal{V(G)}\}$. If
$(v_i,v_j)\in\mathcal{E(G)}$, $v_i$ is called the parent vertex of
$v_j$ and $v_j$ is called the child vertex of $v_i$. The set of {\it
neighbors} of vertex $v_i$ is denoted by $\mathcal{N(G},v_i)=\{v_j:
(v_j,v_i)\in\mathcal{E(G)}, j\not=i\}$. The associated index set is
denoted by $\mathcal{N(G,}i)=\{j: v_j\in\mathcal{N(G,}v_i)\}$.  A
{\it subgraph} $\mathcal{G}_s$ of directed graph $\mathcal{G}$ is a
directed graph such that the vertex set $\mathcal{V}(\mathcal{G}_s)
\subset \mathcal{V}(\mathcal{G})$ and the edge set
$\mathcal{E}(\mathcal{G}_s) \subset \mathcal{E}(\mathcal{G})$.  If
$\mathcal{V}(\mathcal{G}_s) = \mathcal{V}(\mathcal{G})$, we call
$\mathcal{G}_s$ a {\it spanning subgraph} of $\mathcal{G}$. For any
$v_i, v_j\in\mathcal{V}(\mathcal{G}_s)$, if $(v_i,v_j)\in
\mathcal{E}(\mathcal{G}_s)$ if and only if $(v_i,v_j)\in
\mathcal{E}(\mathcal{G})$, $\mathcal{G}_s$ is called an {\it induced
subgraph}.  In this case, $\mathcal{G}_s$ is also said to be induced
by $\mathcal{V}(\mathcal{G}_s)$. A {\it path} in directed graph
$\mathcal{G}$ is a finite sequence $v_{i_1}, \dots, v_{i_j}$ of
vertices such that $(v_{i_k}, v_{i_{k+1}})\in \mathcal{V(G)}$ for
$k=1, \dots, j-1$.     A {\it directed tree} is a directed graph,
where every vertex, except one special vertex without any parent,
has exactly one parent, and the special vertex, called the {\it root
vertex}, can be connected to any other vertices through paths. A
{\it spanning tree} of $\mathcal{G}$ is a directed tree that is a
spanning subgraph of $\mathcal{G}$. We say that a directed graph has
or contains a spanning tree if a subset of the edges forms a
spanning tree. Directed graph $\mathcal{G}$ is {\it strongly
connected} if between every pair of distinct vertices $v_i, v_j$ in
$\mathcal{G}$, there is a path  that begins at $v_i$ and ends at
$v_j$ (that is, from $v_i$ to $v_j$).  A {\it strongly connected
component} of a directed graph is an induced subgraph that is
maximal, subject to being strongly connected. Since any subgraph
consisting of one vertex is strongly connected, it follows that each
vertex lies in a strongly connected component, and therefore the
strongly connected components of a given directed graph  partition
its vertices.

Matrix $\A$ is called a nonnegative matrix, denoted by $\A\geq 0$,
if all its entries are nonnegative, and is called a positive matrix,
denoted by $\A> 0$, if all its entries are positive.
 A {\it weighted
directed graph} $\mathcal{G}(\A)$ is a directed graph $\mathcal{G}$
plus a nonnegative {\it weight matrix} $\A=[a_{ij}]\in \R^{n\times
n}$ such that $(v_i, v_j)\in \mathcal{E(G)} \iff a_{ji}>0$. And
$a_{ji}$ is called the {\it weight} of edge $(v_i, v_j)$. Moreover,
  if $\A^T=\A$, then
$\mathcal{G}(\A)$ is also called an {\it undirected} graph. In this
case, $(v_i, v_j)\in\G(\A)\iff (v_j,v_i)\in\G(\A)$, and
$\mathcal{G}(\A)$ having a spanning tree is equivalent to
$\mathcal{G}(\A)$ being strongly connected. Those strongly connected
undirected graphs are usually called to be {\it connected}. In this
paper, the induced subgraph of weighted directed graph $\G(\A)$ is
also a
 weighted directed graph which inherits its weights in $\G(\A)$.

{\it Notations:} Let $\yi=[1,1,\dots,1]^T$ with compatible
dimensions, $\n=\{1,2,\cdots,n\}$, $\spans(\yi)=\{\bxi\in\R^n:
\bxi=r\yi, r\in\R\}$, and $\rho(\A)$ denote the spectral radius of a
square matrix $\A$. For any $1\leq p\leq \infty$, $\|\cdot\|_{p}$
denotes the $l_p$-norm on $\R^n$ and $\|\cdot\|_{ip}$ is its induced
norm on $\R^{n\times n}$. If
$\w=[\omega_1,\omega_2,\dots,\omega_n]^T\in\R^n$, then $\diag(\w)$
is the diagonal matrix with the $(i,i)$ entry being $\omega_i$.

\begin{yl}[\cite{R. Olfati-Saber and R. M.  Murray
1,W. Ren and R. W. Beard,xfisic2006}]\label{yl2} Let
$L(\A)=[l_{ij}]\in\R^{n\times n}$ denote the {\it graph Laplacian}
of $\mathcal{G}(\A)$, which  is defined by
\[
l_{ij}=\left\{
         \begin{array}{ll}
           \sum_{k=1,k\not=i}^n a_{ik}, & j=i \\
           -a_{ij}, & j\not=i \\
         \end{array}
       \right..
\]
Then

(i) $0$ is an eigenvalue of $L(\A)$  and $\yi$ is the associated
eigenvector;

(ii) If $\mathcal{G}(\A)$ has a spanning tree, then eigenvalue $0$
is algebraically simple and all other eigenvalues are with positive
real parts;

(iii) If $\mathcal{G}(\A)$ is strongly connected, then there exists
a positive column vector $\w\in\R^{ n}$ such that $\w^T L(\A)=0$;

If $\mathcal{G}(\A)$ is undirected, namely, $\A^T=\A$, and
connected, then $L(\A)$ has the following properties:

(iv) $\bxi^TL(\A)\bxi=\frac{1}{2}\sum_{i,j=1}^n
a_{ij}(\xi_j-\xi_i)^2$ for any
$\bxi=[\xi_1,\xi_2,\dots,\xi_n]^T\in\R^n$, and therefore $L(\A)$ is
semi-positive definite, which implies that all eigenvalues of
$L(\A)$ are nonnegative real numbers;

(v) The second smallest eigenvalue of $L(\A)$, which is denoted by
$\lambda_2(L(\A))$ and called the {\it algebraic connectivity} of
$\mathcal{G}(\A)$, is larger than zero;

(vi) The algebraic connectivity of $\mathcal{G}(\A)$ is equal to
$\min_{\bxi\not=0, \yi^T\bxi=0}\frac{\bxi^TL(\A)\bxi}{\bxi^T\bxi}$,
and therefore, if $\yi^T\bxi =0$, then
\[\bxi^TL(\A)\bxi\geq \lambda_2(L(\A))\bxi^T\bxi.\]
\end{yl}

\begin{proof}
Statement (i) follows directly from the definition of  $L(\A)$. (ii)
is a corollary of Lemma 3.3 or Lemma 3.11 in \cite{W. Ren and R. W.
Beard}. The last three statements appeared in \cite{R. Olfati-Saber
and R. M. Murray 1} (see equations (15) (17) and Theorem 1 in
\cite{R. Olfati-Saber and R. M.  Murray 1}).

Here, we provide the proof of (iii) for completeness, which can be
extended to more general case, see \cite{xfisic2006}.

Let $d'=\max_i l_{ii}$, $\varepsilon=\max_ia_{ii}+1$ and let
$d=d'+\varepsilon$. Then $-L(\A)=-d \bI+(-L(\A)+ d \bI)$, where
$\bI$ is the identity matrix with compatible dimensions. By direct
observation, $-L(\A)+ d \bI\geq \A$ (component-wise) is a
nonnegative matrix with positive diagonal entries. Since $\G(\A)$ is
strongly connected, $\G(-L(\A)+d\bI)$ is also strongly connected. By
Lemma \ref{ylHorn6224} (in the Appendix), $-L(\A)+d\bI$ is
irreducible, equivalently, $-L(\A)^T+d\bI$ is also irreducible.

By Ger\v{s}gorin Disk Theorem (in the Appendix), all the eigenvalues
of $-L(\A)$ are located in the following region:
\[
\bigcup_{i=1}^{n}\{c\in\mathbb{C}: |c+l_{ii}|\leq\sum_{j=1\atop
j\not=i}^n|l_{ij}|=l_{ii}\}\subset\{c\in\mathbb{C}:|c+d'|\leq d'\}
\]
\begin{figure}[htpb]\centering
      \includegraphics[scale=0.38]{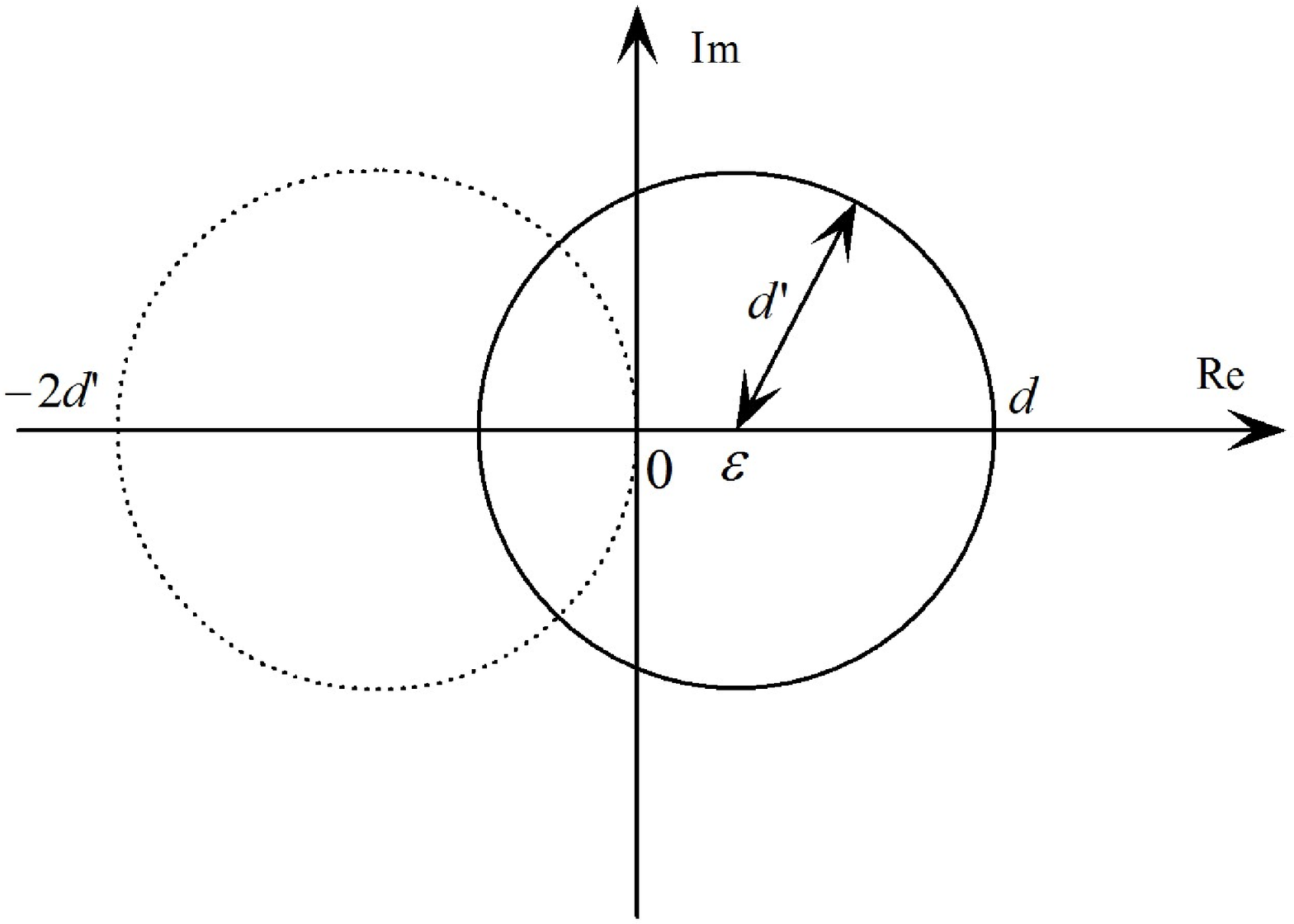}\\
\small The eigenvalues of $-L(\A)+  d\bI$ are located in the right
circle.
    \caption{Eigenvalues of $-L(\A)+  d\bI$ }
      \label{eigenvaluelacation}
\end{figure}

Therefore, the eigenvalue of $-L(\A)+d\bI$ are located in the region
$\{c\in\mathbb{C}: |c-\varepsilon|\leq d'\}$, see Fig.
\ref{eigenvaluelacation}. Consequently, its spectral radius,
$\rho(-L(\A)+d\bI)$, is not larger than $d$. From
$(-L(\A)+d\bI)\yi=(0+d)\yi=d\yi$, we have that $d$ is an eigenvalue
of $-L(\A)+d\bI$, and thus $\rho(-L(\A)+d\bI)=d$.  By
Perron-Frobenius Theorem (in the Appendix), there exists a positive
column vector $\w$ such that $(-L(\A)^T+d \bI)\w=d\w$. Therefore,
$\w^T L(\A)=\w^T (d\bI) -d\w^T=0$.
\end{proof}

\begin{tl}\label{tlwL}
Suppose $\G(\A)$ is strongly connected, and let $\w>0$ such that
$\w^T L(\A)=0$. Then $\diag(\w)L(\A)+L(\A)^T\diag(\w)$ is the graph
Laplacian of the undirected weighted graph
$\G(\diag(\w)\A+\A^T\diag(\w))$. And therefore it is semi-positive
definite,  $0$ is its algebraically simple eigenvalue and $\yi$ is
the associated eigenvector.
\end{tl}
\begin{proof}
Since the definition of graph Laplacian $L(\A)$ is independent of
the diagonal entries of $\A$, without loss of generality, we assume
that the diagonal entries of $\A$ are all zeros (in the subsequent
sections, $\A$ will be assumed to be with zero diagonal entries).
Then by the definition of graph Laplacian, we have that
\begin{equation}\label{laplacianLyi1}
    L(\A)=\diag(\A\yi)-\A.
\end{equation}
Then\begin{align}
    \diag(\w)L(\A)=&\diag(\w)(\diag(\A\yi)-\A)\nonumber\\
    =&\diag(\diag(\w)\A\yi)-\diag(\w)\A\nonumber\\
    =&L(\diag(\w)\A)\mbox{ (by \eqref{laplacianLyi1}) }\label{laplacianLyi2}.
\end{align}
Since $\w^TL(\A)=0$, $\yi^T\diag(\w)L(\A)=0$. Therefore,
\begin{align*}
    0=&\yi^T\diag(\w)L(\A)=\yi^TL(\diag(\w)\A)\\
    =&\yi^T\diag(\diag(\w)\A\yi)-\yi^T(\diag(\w)\A)\\
    =&\yi^T\A^T\diag(\w)-\yi^T(\diag(\w)\A).
\end{align*}
Thus,
\begin{equation}\label{laplacianLyi3}
    \diag(\w)\A\yi=\A^T\diag(\w)\yi.
\end{equation}

Therefore,
\begin{align*}
L(\A)^T\diag(\w)=&L(\diag(\w)\A)^T\mbox{ (by \eqref{laplacianLyi2}) }=\diag(\diag(\w)\A\yi)-(\diag(\w)\A)^T\\
=&\diag(\A^T\diag(\w)\yi)-\A^T\diag(\w) \mbox{ (by \eqref{laplacianLyi3})}\\
 =&L(\A^T\diag(\w)) \mbox{  (by \eqref{laplacianLyi1}) }.
\end{align*}

Moreover, by \eqref{laplacianLyi1}, we can see that for any
nonnegative matrices $\A_1$ and $\A_2$,
$L(\A_1)+L(\A_2)=L(\A_1+\A_2)$. Therefore,
\begin{align*}
    \diag(\w)L(\A)+&L(\A)^T\diag(\w)=L(\diag(\w)\A)+L(\A^T\diag(\w))\\
    &=L(\diag(\w)\A+\A^T\diag(\w)).
\end{align*}
\end{proof}

\begin{tl}\label{tlLb}
Let $\bb=[b_1,b_2,\dots,b_n]^T\geq 0$, $\bb\not=0$, and let $\G(\A)$
be undirected and connected. Then $L(\A)+\diag(\bb)$ is positive
definite.
\end{tl}
\begin{proof}
By Lemma \ref{yl2}(iv), $L(\A)+\diag(\bb)$ is semi-positive
definite. If there exists a vector $\bxi\in\R^n$ such that
$\bxi^T(L(\A)+\diag(\bb))\bxi=0$, then by lemma \ref{yl2}, we have

(i) $\bxi\in\spans(\yi)$;

(ii) $\bxi^T\diag(\bb)\bxi=0$, which implies that some entries of
$\bxi$ are  zeros.

Hence, $\bxi=0$ and $L(\A)+\diag(\bb)$ is positive definite.
\end{proof}

\section{Problem Formulation}
 The distributed dynamic system studied in this paper
consists of $n$ autonomous agents, e.g. particles or robots, labeled
$1$ through $n$. All these agents share a common state space $\R$.
The state of agent $i$ is denoted by $x_i$, $i=1,2,\dots,n$ and let
$\x=[x_1, x_2,\dots, x_n]^T$.

Suppose that agent $i$, $i\in\n$, is with the following dynamics
\begin{equation}\label{xt}
    \dot{x}_i(t)=u_i(t),
\end{equation}
where $u_i(t)$   is the protocol to be designed.

In this multi-agent system, each agent can communicate with some
other agents which are defined as its neighbors. Protocol $u_i$ is a
state feedback, which is designed based on the state information
received by agent $i$ from its neighbors. We use the weighted
directed graph $\G(\A)$ to represent the communication topology,
where $\A=[a_{ij}]$ is a given $n\times n$ nonnegative matrix, with
diagonal entries all being zeros. Vertex $v_i$ represents agent $i$.
Edge $(v_i, v_j)\in\bN(\G(\A))$ represents an available information
channel from agent $i$ to agent $j$. The {\it neighbors} of agent
$i$ are those agents whose information can be received directly by
agent $i$. Thus, the set of neighbors of agent $i$ just corresponds
to the set $\bN(\G(\A), v_i)$. The {\it leaders } of this system are
the agents,  the vertices corresponding to which are the roots of
some spanning trees of $\G(\A)$, i.e., the vertices that can be
connected to any other vertices through paths \cite{xfjmaa}. The
agents that are not leaders are called followers. Given certain
agents, the {\it local communication topology } of them is the
subgraph of $\G(\A)$, which is induced by the vertices corresponding
to those agents. If $\A$ is time-dependent, then the underlying
topology is said to be switching. In order to reflect the dependency
on time, we use $\A(t)$ instead of $\A$ in this case. The switching
topology will be
 discussed in detail in Section IV.D. Unless it is explicitly specified,
$\A$ is assumed to be time-invariant.

Given protocol $u_i$, $i\in\n$, we say that $u_i$ or this
multi-agent system solves a consensus or agreement problem if for
any given initial states and any $j,k\in\n$, $|x_j(t)-x_k(t)|\to 0$,
as $t\to\infty$, and we say that it solves a finite-time consensus
problem if it solves a consensus problem, and given any initial
states, there exist a time $t^*$ and a real number $\kappa$ such
that $x_j(t)=\kappa$ for $t\geq t^*$ and  all $j\in\n$. If the final
consensus state is a function of the initial states, namely,
$x_j(t)\to \chi(\x(0))$ for all $j\in\n$ as $t\to\infty$, where
$\chi:\R^n\to\R$ is a function, we say that it solves the
$\chi$-consensus problem \cite{R. Olfati-Saber and R. M.  Murray 1}.
Specially, if $\chi(x(0))=\frac{\sum_{i=1}^nx_i(0)}{n}$, the system
is said to solve the average-consensus problem.

We are now in a position to present two consensus protocols, which
will be shown to  solve finite-time consensus problems:

(i) For $i\in\n$, if $\mathcal{N}(\G(\A), v_i)=\phi$, then $u_i=0$,
else
\begin{equation}\label{pro1}
    u_i=\sign\left(\sum_{j\in
    \bN(\G(\A),i)}a_{ij}(x_j-x_i)\right)\left|\sum_{j\in
    \bN(\G(\A),i)}a_{ij}(x_j-x_i)\right|^{\alpha_i};
\end{equation}

 (ii) For $i\in\n$, if $\mathcal{N}(\G(\A), v_i)=\phi$, then $u_i=0$, else
\begin{equation}\label{pro2}
    u_i=\sum_{j\in
    \bN(\G(\A),i)}a_{ij}\sign(x_j-x_i)|x_j-x_i|^{\alpha_{ij}},
\end{equation}
 where $0<\alpha_i, \alpha_{ij}<1$, $|\cdot|$ is   the absolute value of real numbers, and $\sign(\cdot)$ is the sign
function defined as
\begin{equation*}
    \sign(r)=\left\{
            \begin{array}{rl}
              1, & r>0 \\
              0, & r=0 \\
              -1, & r<0 \\
            \end{array}.
          \right.
\end{equation*}

For simplicity, $\sign(r)|r|^{\alpha}$ is denoted by
$\sgn{r}{\alpha}$ for any $\alpha>0$.

{\it Remark:} $\sgn{r}{\alpha}$, $\alpha>0$, is a continuous
function on $\R$, which leads to the continuity of protocols
\eqref{pro1} and \eqref{pro2}. That will be discussed in Property
\ref{xz1}. If we set $\alpha_i=1,i\in\n$ and $\alpha_{ij}=1,
i,j\in\n$ in the above protocols \eqref{pro1} and \eqref{pro2}
respectively, then they will become the typical linear consensus
protocol studied in \cite{R. Olfati-Saber and R. M. Murray 1} and
\cite{W. Ren and R. W. Beard}, and they solve a consensus problem
asymptotically provided that $\mathcal{G}(\A)$ has a spanning tree.
If we set $\alpha_i=0$ and $\alpha_{ij}=0$  in \eqref{pro1} and
\eqref{pro2}, they will become discontinuous.  The discontinuous
case of protocol \eqref{pro1} was studied by Cort\'{e}s in \cite{J.
Cortes}, where $\G(\A)$ is assumed to be a connected undirected
graph and $\A$ is also a $0-1$ matrix. In the next section, we will
investigate the mathematical conditions that guarantee protocol
\eqref{pro1} or protocol \eqref{pro2} to solve a finite-time
consensus problem.

\subsection{Basic Properties and Lemmas}

\begin{xz}\label{xz1}
Protocols \eqref{pro1} and \eqref{pro2} are continuous with respect
to state variables $x_1,x_2, \dots, x_n$. Moreover, under either of
those two protocols, there exists at least one solution of
differential equations \eqref{xt} on $[0, \infty)$ for any initial
state $\x(0)$. Furthermore, $\max_i x_i(t)$ is non-increasing and
$\min_i x_i(t)$ is non-decreasing. Hence, $\|\x(t)\|_{\infty}$ is
also non-increasing and $\|\x(t)\|_{\infty}\leq \|\x(0)\|_{\infty}$
for all $t\geq 0$.
\end{xz}
\begin{proof}
Consider function $f(r)=\sgn{r}{\alpha}$, where $0<\alpha<1$. If
$r>0$, then $f(r)=r^{\alpha}$ is continuous, and $\lim_{r\to 0^{+}}
r^{\alpha}=0$. If $r<0$, then $f(r)=-(-r)^{\alpha}$ is also
continuous,  and $\lim_{r\to 0^{-}} -(-r)^{\alpha}=0$. From
$0^{\alpha}=0$, we have that $f(r)$ is continuous on $\R$, and
therefore, protocol \eqref{pro1} and \eqref{pro2} are continuous
with respect to state variables $x_1,x_2, \dots, x_n$.

Let $E$ be the set $\{(t,\bxi): \|\bxi\|_{\infty}<
\|\x(0)\|_{\infty}+1, -r<t<r, r>0\}$. Then $E$ is open. Consider the
solution when the domain of $u_i$ is restricted in the open set $E$.
By Peano's Existence Theorem (in the Appendix), there exists at
least one solution of differential equations \eqref{xt} on $[0,t_+]$
for some $t_+>0$. We extend the obtained solution over a maximal
interval of existence $(\omega_-, \omega_+)$. At time $t$,
$t\in[0,\omega_+)$, if $x_j(t)=\max_i x_i(t)$, then
$\dot{x}_j(t)\leq 0$, and if $x_j(t)=\min_i x_i(t)$, then
$\dot{x}_j(t)\geq 0$. Therefore $\max_i x_i(t)$ is non-increasing
and  $\min_i x_i(t)$ is non-decreasing. With the same arguments,
$\|\x(t)\|_{\infty}$ is also non-increasing for $0\leq t<\omega_+$.
By Extension Theorem (in the appendix), $(t,\x(t))$ tends to the
boundary of $E$ as $t\to\omega_+$, which only occurs when
$\omega_+=r$. Because of the arbitrariness of $r$, there exists at
least one solution of differential equation \eqref{xt} on $[0,
\infty)$ for any initial state $\x(0)$.
\end{proof}

{\it Remark:} In fact,  any solution of the system can be extended
over $[0,\infty)$ because it is  bounded by $\|\x(0)\|_{\infty}$.

Moreover, one notices that \eqref{pro1} and \eqref{pro2} are not
Lipschitz at some points. As all solutions reach subspace
$\spans(\yi)$ in finite time, there is nonuniqueness of solutions in
backwards time. This, of course, violates the uniqueness condition
for solutions of Lipschitz differential equations.

The following property shows that the equilibrium point set of the
considered differential equations $\dot{x}_i=u_i, i\in\n$, is the
set of all consensus states.

\begin{xz}\label{xz2}
With protocol \eqref{pro1} or \eqref{pro2}, the equilibrium point
set of the differential equations $\dot{x}_i=u_i, i\in\n$, is
$\spans(\yi)$, provided that $\G(\A)$ has a spanning tree.
\end{xz}
\begin{proof}
Since $\G(\A)$ has a spanning tree, there exists at most one agent,
whose neighbor set is empty. Let $\bxi=[\xi_1, \xi_2, \dots,
\xi_n]^T$ satisfy the equations $u_i|_{\x=\bxi}=0$, $i\in\n$. We
    prove $\bxi\in\spans(\yi)$ by contradiction. Assume that $\min_i\xi_i\not=\max_i\xi_i$. Let
    $\mathcal{M}=\{j:\xi_j=\min_i\xi_i\}$ and let
    $\mathcal{H}=\{j:\xi_j=\max_i\xi_i\}$. Then
    $\mathcal{M}\cap\mathcal{H}=\phi$. Suppose that the root of the spanning
    tree is vertex $v_j$. Then $j\not\in\mathcal{M}$ or
    $j\not\in\mathcal{H}$ or both hold. Without loss of generality,
    assume that $j\not\in\mathcal{M}$. For any $k\in\mathcal{M}$,
    there exists a path $v_j=v_{i_1}, v_{i_2}, \dots, v_{i_s}=v_k$
    from $v_j$ to $v_k$. Let $l$ be the number, such that
    $i_l\in\mathcal{M}$ but $i_{l-1}\not\in\mathcal{M}$. Since $\xi_{i_{l-1}}>\min_i\xi_i$, $i_{l-1}\in\bN(\G,
    i_l)$ and $\xi_{i_l}=\min_i\xi_i$, we have that $u_{i_l}>0$, which is a
    contradiction. Thus, $\min_i\xi_i=\max_i\xi_i$, namely,
    $\bxi\in\spans(\yi)$.

    To prove the converse, suppose that $\bxi\in\spans(\yi)$. Then
    all $u_i|_{\x=\bxi}$, $i\in\n$, are all zeros. Therefore, $\bxi$ belongs to
    the equilibrium point set.
\end{proof}

In order to establish our main results, we need the following
Lemmas.
\begin{yl}\label{yl1}
Let $\xi_1, \xi_2, \dots, \xi_n \geq 0$ and let $0<p\leq 1$. Then
\[\left( \sum_{i=1}^n \xi_i\right)^p\leq \sum_{i=1}^n \xi_i^p\leq n^{1-p}\left( \sum_{i=1}^n \xi_i\right)^p.\]
\end{yl}
\begin{proof}
Obviously
\[ \sum_{i=1}^n \xi_i=0 \iff \sum_{i=1}^n \xi_i^p=0.\]

Let $\mathcal{U}=\{\zt=[\zeta_1,\zeta_2,\dots,\zeta_n]^T\in\R^n:
\sum_{i=1}^n\zeta_i =1\mbox{\ and \ }\zt\geq 0\}$. If $\sum_{i=1}^n
\xi_i\not=0$,
\begin{align*}
    \frac{\sum_{i=1}^n \xi_i^p}{\left( \sum_{i=1}^n
    \xi_i\right)^p}&=\sum_{i=1}^n \left(\frac{\xi_i}{ \sum_{i=1}^n
    \xi_i}\right)^p \geq \inf_{\zt\in \mathcal{U}}\sum_{i=1}^n \zeta_i^p,
\end{align*}
where the last inequality follows from that
$\frac{\bxi}{\yi^T\bxi}\in\mathcal{U}$.

For any $\zt \in\mathcal{U}$,  $\sum_{i=1}^n \zeta_i^p\not=0$. Since
$\sum_{i=1}^n \zeta_i^p$ is continuous and $\mathcal{U}$ is a
bounded closed set, $\inf_{\zt\in \mathcal{U}}\sum_{i=1}^n
\zeta_i^p$ exists and is larger than $0$. It can be calculated
directly. Precisely,
\[\inf_{\zt\in \mathcal{U}}\sum_{i=1}^n \zeta_i^p=\min\{n^{1-p}, 1\}=1.\] Thus, the left inequality holds.

With the same arguments,
\[\frac{\sum_{i=1}^n \xi_i^p}{\left( \sum_{i=1}^n
    \xi_i\right)^p}\leq \sup_{\zt\in \mathcal{U}}\sum_{i=1}^n \zeta_i^p=\max\{n^{1-p}, 1\}=n^{1-p}.\]
Therefore, the right inequality also holds.
\end{proof}

\begin{yl}[cf. \cite{S. J. Bhat ACC}, Theorem 1]\label{ylfinitetimeV}
Suppose that function $V(t):[0, \infty)\to [0, \infty)$ is
differentiable (the derivative of $V(t)$ at $0$ is in fact its right
derivative), such that
\[\frac{dV(t)}{dt}\leq -K V(t)^{\alpha},\]
where $K>0$ and $0<{\alpha}<1$. Then $V(t)$ will reach zero at
finite time $t^*\leq\frac{V(0)^{1-{\alpha}}}{K(1-{\alpha})}$, and
$V(t)=0$ for all $t\geq t^*$.
\end{yl}
\begin{proof}
Let $f(t)$ satisfy differential equation
\[\frac{d f(t)}{dt}=-K f(t)^{\alpha}.\]
Given initial value $f(0)=V(0)>0$, its unique solution is
\[
f(t)=\left\{
  \begin{array}{ll}
    \left(-K(1-\alpha)t+f(0)^{1-{\alpha}}\right)^{\frac{1}{1-{\alpha}}},& t< \frac{V(0)^{1-{\alpha}}}{K(1-{\alpha})},\\
    0, & t\geq \frac{V(0)^{1-{\alpha}}}{K(1-{\alpha})} \\
  \end{array}
\right.
\]

Since $V(0)=f(0)$, by Comparison Principle of differential equations
(in the Appendix), $V(t)\leq f(t)$, $t\geq 0$. Hence, $V(t)$ will
reach zero in time $\frac{V(0)^{1-{\alpha}}}{K(1-{\alpha})}$. Since
$V(t)\geq 0$ and $\frac{d V(t)}{dt}\leq 0$, $V(t)$ cannot leave zero
once it reaches it.
\end{proof}

\section{Main Results}
\subsection{Networks Under Protocol \eqref{pro1}}

In this subsection,  protocol \eqref{pro1} is studied. The
communication topology $\mathcal{G}(\A)$ is supposed to have a
spanning tree, which is the least requirement when the topology is
time-invariant. Otherwise, if $\mathcal{G}(\A)$ does not have a
spanning tree, then there exist at least two groups of agents,
between which there does not exist information exchange, and
therefore it is impossible for the system to solve a consensus
problem through distributive protocols. Furthermore, we hope that
the state of the system reaches $\spans{(\yi)}$ in a finite time.

Now, we present our first main result.

\begin{dl}\label{dl1}
If the communication topology $\mathcal{G}(\A)$ has a spanning tree,
then system \eqref{xt} solves a finite-time consensus problem   when
  protocol \eqref{pro1} is applied.
\end{dl}
\begin{proof}
This theorem is proved through the following three steps.

{\it Step 1:} Suppose that the communication topology
$\mathcal{G}(\A)$ is strongly connected.

By Lemma \ref{yl2}, there exists a vector
$\w=[\omega_1,\omega_2,\dots,\omega_n]^T\in \R^{n}$ such that $\w>0$
and $\w^T L(\A)=0$. Let $y_i=\sum_{j=1}^n a_{ij}(x_j-x_i), i\in\n$,
and let $\y=[y_1,y_2,\dots,y_n]^T$. Then
\[\y=-L(\A)\x,\]
and
\[\dot{x}_i=\sgn{y_i}{\alpha_i}.\]
Since $\w^T\y(t)=-\w^T L(\A)\x(t)=-0\x(t)=0,
$
 $\y(t)\perp \w$. Consider nonnegative function
\[V_1(t)=\sum_{i=1}^n
\frac{\omega_i}{\alpha_i+1}|y_i|^{\alpha_i+1},\] which will be
proved to be a valid Lyapunov function for the disagreement of
agents' states.

 If $y_i<0$,
$|y_i|^{\alpha_i+1}=(-y_i)^{\alpha_i+1}$ and $\frac{d
(-y_i)^{\alpha_i+1}
}{dy_i}=-(\alpha_i+1)(-y_i)^{\alpha_i}=(\alpha_i+1)\sgn{y_i}{\alpha_i}$.
If $y_i> 0$, $|y_i|^{\alpha_i+1}=y_i^{\alpha_i+1}$ and $\frac{d
y_i^{\alpha_i+1}
}{dy_i}=(\alpha_i+1)y_i^{\alpha_i}=(\alpha_i+1)\sgn{y_i}{\alpha_i}$.
Furthermore, the left and right derivatives of $|y_i|^{\alpha_i+1}$
at $0$ are all $0$. Consequently, for any $y_i$,
\[\frac{d |y_i|^{\alpha_i+1}
}{dy_i}=(\alpha_i+1)\sgn{y_i}{\alpha_i}.\] Thus,
\begin{align*}
    \frac{d V_1(t)}{dt}&=\sum_{i=1}^n \omega_i \sgn{y_i}{\alpha_i}\frac{d
y_i}{dt}\\
&=\sum_{i=1}^n\left(\omega_i \sgn{y_i}{\alpha_i}\sum_{j=1}^n
a_{ij}\left(\sgn{y_j}{\alpha_j} - \sgn{y_i}{\alpha_i}\right)\right).
\end{align*}

Let $\al=[\alpha_1,\alpha_2,\dots,\alpha_n]^T$, let $\sgn{\y}{\al} =
[\sgn{y_1}{\alpha_1}, \sgn{y_2}{\alpha_2},\dots,
\sgn{y_n}{\alpha_n}]^T$ and let $\sgn{\y^T}{\al}=(\sgn{\y}{\al})^T$.
Then
\[\frac{V_1(t)}{dt}=-\sgn{\y^T}{\al}\diag(\w)L(\A)\sgn{\y}{\al}.\]

Let $\alpha_0=\max_{i}\alpha_i$ and let
$\B=\frac{1}{2}(\diag(\w)L(\A)+L(\A)^T\diag(\w))$. Suppose that
$V_1(t)\not=0$. Then $\y\not=0$, and
\begin{align}
\frac{d V_1(t)}{dt}&=-\frac{1}{2}\sgn{\y^T}{\al}(\diag(\w)L(\A)+L(\A)^T\diag(\w))\sgn{\y}{\al}\nonumber\\
&=-\frac{\sgn{\y^T}{\al}\B\sgn{\y}{\al}}{\sgn{\y^T}{\al}\sgn{\y}{\al}}
\frac{\sgn{\y^T}{\al}\sgn{\y}{\al}}{V_1(t)^{\frac{2\alpha_0}{1+\alpha_0}}}V_1(t)^{\frac{2\alpha_0}{1+\alpha_0}}\label{eqdl1Vfrac}
\end{align}
Since $0<\alpha_0<1$, $0<\frac{2\alpha_0}{1+\alpha_0}<1$. In order
to apply Lemma \ref{ylfinitetimeV}, we need to find lower bounds of
the first two qualities in the right side of equality
\eqref{eqdl1Vfrac}.

Let $\mathcal{U}=\{\bxi: \mbox{nonzero terms of }\xi_1,
\xi_2,\dots,\xi_n$ are not with the same sign $ \}$ and let
$\mathcal{U}^0= \mathcal{U}\cap\{\bxi:\bxi^T\bxi=1\}$. Obviously,
$\mathcal{U}^0$ is a bounded closed set, namely, a compact set. From
$\y\perp\w$ and $\w>0$, we have $\y\in\mathcal{U}$, and from that
function $\sgn{r}{\alpha}$ conserves the sign of $r$, we have
$\sgn{\y}{\al}\in\mathcal{U}$. Let $\bxi\in\mathcal{U}^0$. Then
$\bxi\not\in\spans(\yi)$ and therefore by Corollary \ref{tlwL},
$\bxi^T\B\bxi>0$. Because the function $\bxi^T\B\bxi$ is continuous
and $\mathcal{U}^0$ is compact,
$\min_{\bxi\in\mathcal{U}^0}\bxi^T\B\bxi$, denoted by $K_1$, exists
and is larger than zero. Hence
\begin{align*}
    \frac{\sgn{\y^T}{\al}\B\sgn{\y}{\al}}{\sgn{\y^T}{\al}\sgn{\y}{\al}}
    =\frac{\sgn{\y^T}{\al}}{\sqrt{\sgn{\y^T}{\al}\sgn{\y}{\al}}}\B \frac{\sgn{\y}{\al}}{\sqrt{\sgn{\y^T}{\al}\sgn{\y}{\al}}}
    \geq K_1.
\end{align*}
The last inequality follows from that
$\frac{\sgn{\y}{\al}}{\sqrt{\sgn{\y^T}{\al}\sgn{\y}{\al}}}\in\mathcal{U}^0$.

Consider the second quality,
\begin{align}
    \frac{\sgn{\y^T}{\al}\sgn{\y}{\al}}{V_1(t)^{\frac{2\alpha_0}{1+\alpha_0}}}&=
    \frac{\sum_{i=1}^n
    |y_i|^{2\alpha_i}}{\left(\sum_{i=1}^n\frac{\omega_i}{1+\alpha_i}|y_i|^{1+\alpha_i}\right)^{\frac{2\alpha_0}{1+\alpha_0}}} \nonumber\\
    &\geq \frac{\sum_{i=n}^n
    |y_i|^{2\alpha_i}}{\sum_{i=1}^n
    \left(\frac{\omega_i}{1+\alpha_i}\right)^{\frac{2\alpha_0}{1+\alpha_0}}|y_i|^{(1+\alpha_i)\frac{2\alpha_0}{1+\alpha_0}}}\mbox{ ( by Lemma \ref{yl1})}\nonumber\\
    &\geq
    \frac{|y_l|^{2\alpha_l}}{\left(\sum_{i=1}^n \left(\frac{\omega_i}{1+\alpha_i}\right)^{\frac{2\alpha_0}{1+\alpha_0}}\right)|y_l|^{(1+\alpha_l)\frac{2\alpha_0}{1+\alpha_0}}}\nonumber\\
    &=\frac{1}{\sum_{i=1}^n\left(\frac{\omega_i}{1+\alpha_i}\right)^{\frac{2\alpha_0}{1+\alpha_0}}}|y_l|^{2\alpha_l-(1+\alpha_l)\frac{2\alpha_0}{1+\alpha_0}},\nonumber
\end{align}
where
$l=\myarg\max_i|y_i|^{(1+\alpha_i)\frac{2\alpha_0}{1+\alpha_0}}$.

Since $0< \alpha_l\leq\alpha_0< 1$,
$\frac{2\alpha_0}{1+\alpha_0}-\frac{2\alpha_l}{1+\alpha_l}=\frac{2\alpha_0-2\alpha_l}{(1+\alpha_0)(1+\alpha_l)}\geq
0$, which implies that
$2\alpha_l-(1+\alpha_l)\frac{2\alpha_0}{1+\alpha_0}\leq 0$.
Therefore, function
$|r|^{2\alpha_l-(1+\alpha_l)\frac{2\alpha_0}{1+\alpha_0}}$ does not
increase as $|r|$ increases. Since $0<|y_l(t)|\leq
\|\y(t)\|_{\infty}=\|-L(\A)\x(t)\|_{\infty}\leq
\|L(\A)\|_{i\infty}\|\x(t)\|_{\infty}\leq
\|L(\A)\|_{i\infty}\|\x(0)\|_{\infty}$,
$|y_l|^{2\alpha_l-(1+\alpha_l)\frac{2\alpha_0}{1+\alpha_0}}\geq
\left(\|L(\A)\|_{i\infty}\|\x(0)\|_{\infty}\right)^{2\alpha_l-(1+\alpha_l)\frac{2\alpha_0}{1+\alpha_0}}$.
 Let
\[K_2=\frac{1}{\sum_{i=1}^n(\frac{\omega_i}{1+\alpha_i})^{\frac{2\alpha_0}{1+\alpha_0}}}\min_{i}\left(\|L(\A)\|_{i\infty}\|\x(0)\|_{\infty}\right)^{2\alpha_i-(1+\alpha_i)\frac{2\alpha_0}{1+\alpha_0}}.\]
Then $K_2>0$ and
\begin{equation}\label{eqdl1inq2}
    \frac{\sgn{\y^T}{\al}\sgn{\y}{\al}}{V_1(t)^{\frac{2\alpha_0}{1+\alpha_0}}}\geq
    K_2.
\end{equation}
Thus by \eqref{eqdl1Vfrac},
\[\frac{d V_1(t)}{dt}\leq - K_1K_2
V_1(t)^{\frac{2\alpha_0}{1+\alpha_0}}, t\geq 0.\] By Lemma
\ref{ylfinitetimeV}, the above differential inequality  gives that
$V_1(t)$ reaches zero in  finite time
$\frac{(1+\alpha_0)V_1(0)^{\frac{1-\alpha_0}{1+\alpha_0}}}{K_1K_2(1-\alpha_0)}$.
If $V_1(t)=0$, then $\y(t)=0$, which implies that ${u}_i=0, i\in\n$,
and thus $\x(t)\in\spans(\yi)$ by Property \ref{xz2}. Therefore the
system solves a finite-time consensus problem.

{\it Step 2:} Next, suppose that there exists only one leader and
the local communication topology among the followers are strongly
connected. By protocol \eqref{pro1}, the state of the leader is
time-invariant. Without loss of generality, suppose that agents
$1,2,\dots, m=n-1$ are the followers and agent $n$ is the leader
(because it is just a matter of relabeling the $n$ agents).
Consequently, $a_{n1}, a_{n2}, \dots, a_{nn}$ are all equal to zero.
For notational simplicity, let $b_1=a_{1n}, b_2=a_{2n}, \dots,
b_m=a_{mn}$, let $\bar{\bb}=[b_1,b_2,\dots,b_m]^T$, let
$\bar{\al}=[\alpha_1,\alpha_2,\dots,\alpha_m]$, and let
$\bar{\A}=[a_{ij}]_{1\leq i,j\leq m}\in\R^{m\times m}$. Since agent
$n$ is the leader, $\bar{\bb}\not=0$. Denote $x_i-x_n$ by
$\bar{x}_i$, $i\in \m$. Then for any $i\in\m$, we have that
\[\dot{\bar{x}}_i=\dot{x}_i=\sgn{\sum_{j=1}^ma_{ij}(\bar{x}_j-\bar{x}_i)-b_i\bar{x}_i}{\alpha_i}.\]

Let $y_i=\sum_{j=1}^m a_{ij}(\bar{x}_j-\bar{x}_i)-b_i\bar{x}_i$,
$i\in\m$, and let $\bar{\y}=[y_1,y_2,\dots ,y_m]^T$. Then
\[\dot{y}_i=\sum_{j=1}^ma_{ij}(\sgn{y_j}{\alpha_j}-\sgn{y_i}{\alpha_i})-b_i\sgn{y_i}{\alpha_i}, i\in\m.\]

Since the local communication topology among the followers is
$\G(\bar{\A})$ that is strongly connected, by Lemma \ref{yl2}, there
exists $\bar\w=[\omega_1,\omega_2,\dots ,\omega_m]^T\in\R^{m}$ such
that $\bar{\w}>0$ and $\bar\w^T L(\bar\A)=0$. Define Lyapunov
function
$V_2(t)=\sum_{i=1}^m\frac{\omega_i}{1+\alpha_i}|y_i|^{1+\alpha_i}$.
Obviously, $V_2(t)=0$ if and only if $\bar{\y}=0$. Notice that
$\diag(\bar\w)\diag(\bar\bb)=\diag([\omega_1b_1,\omega_2b_2,\dots,\omega_mb_m]^T)$.
Then
$\diag(\bar\w)(L(\bar\A)+\diag(\bar\bb))+(L(\bar\A)^T+\diag(\bar\bb))\diag(\bar\w)$
is positive definite by Corollaries \ref{tlwL} and \ref{tlLb}, and
thus $L(\bar\A)+\diag(\bar\bb)$ is non-degenerate. And from
$\bar\y=-(L(\bar\A)+\diag(\bar\bb))[\bar x_1,\bar x_2,\dots, \bar
x_m]^T$, we have that $[\bar x_1,\bar x_2,\dots, \bar x_m]^T=0$ if
and only if $\bar\y=0$. Therefore, $V_2(t)=0$ if and only
$\x(t)\in\spans(\yi)$.

Differentiate it with respect to time.
\begin{align*}
\frac{dV_2(t)}{dt}&=-\sgn{\bar{\y}^T}{\bar\al}\diag(\bar{\w})L(\bar{\A})\sgn{\bar{\y}}{\bar\al}-\sgn{\bar{\y}^T}{\bar\al}\diag(\bar\w)\diag(\bar\bb)\sgn{\bar{\y}}{\bar\al}\\
&=-\sgn{\bar{\y}^T}{\bar\al}\left(\frac{1}{2}\left(\diag(\bar{\w})L(\bar{\A})+L(\bar{\A})^T\diag(\bar{\w})\right)+\diag(\bar\w)\diag(\bar\bb)\right)\sgn{\bar{\y}}{\bar\al}
\end{align*}
Let
$\bar{\B}=\frac{1}{2}\left(\diag(\bar{\w})L(\bar{\A})+L(\bar{\A})^T\diag(\bar{\w})\right)+\diag(\bar\w)\diag(\bar\bb)$
and $\bar{\alpha}_0=\max_{i\in\m}\alpha_i$. Suppose that
$V_2(t)\not=0$.
\begin{align*}
    \frac{dV_2(t)}{dt}&\leq -\frac{\sgn{\bar{\y}^T}{\bar\al}\bar\B\sgn{\bar{\y}}{\bar\al}}{\sgn{\bar{\y}^T}{\bar\al}\sgn{\bar{\y}}{\bar\al}}
    \frac{\sgn{\bar{\y}^T}{\bar\al}\sgn{\bar{\y}}{\bar\al}}{V_2(t)^{\frac{2\bar{\alpha}_0}{1+\bar{\alpha}_0}}}V_2(t)^{\frac{2\bar{\alpha}_0}{1+\bar{\alpha}_0}}.
\end{align*}

By Corollaries \ref{tlwL} and \ref{tlLb}, $\bar{\B}$ is real
symmetric and positive definite. Let the smallest eigenvalue of
$\bar{\B}$ be $\lambda_1(\bar{\B})$. Then $\lambda_1(\bar{\B})>0$,
and for any $\bxi\in \R^m$, $\bxi^T\bar\B\bxi\geq
\lambda_1(\bar{\B})\bxi^T\bxi$. Therefore,
\[\frac{\sgn{\bar{\y}^T}{\bar\al}\bar\B\sgn{\bar{\y}}{\bar\al}}{\sgn{\bar{\y}^T}{\bar\al}\sgn{\bar{\y}}{\bar\al}}\geq \lambda_1(\bar\B).\]

With the same arguments as in the proof of inequality
\eqref{eqdl1inq2}, there exists
\[K_3=\frac{1}{\sum_{i=1}^m\left(\frac{\omega_i}{1+\alpha_i}\right)^{\frac{2\bar{\alpha}_0}{1+\bar{\alpha}_0}}}\min_{i\in\m}
\left(\|L(\A)\|_{i\infty}\|\x(0)\|_{\infty}\right)^{2\alpha_i-(1+\alpha_i)\frac{2\bar{\alpha}_0}{1+\bar{\alpha}_0}}>0\]
such that
\[
\frac{\sgn{\bar{\y}^T}{\bar\al}\sgn{\bar{\y}}{\bar\al}}{V_2(t)^{\frac{2\bar{\alpha}_0}{1+\bar{\alpha}_0}}}\geq
K_3.
\] Therefore,
\[\frac{d V_2(t)}{dt}\leq -\lambda_1(\bar\B)K_3V_2(t)^{\frac{2\bar{\alpha}_0}{1+\bar{\alpha}_0}}.\]
By Lemma \ref{ylfinitetimeV}, in this case, this system solves a
finite-time consensus problem.

\begin{figure}[htpb]
      \includegraphics[scale=0.8]{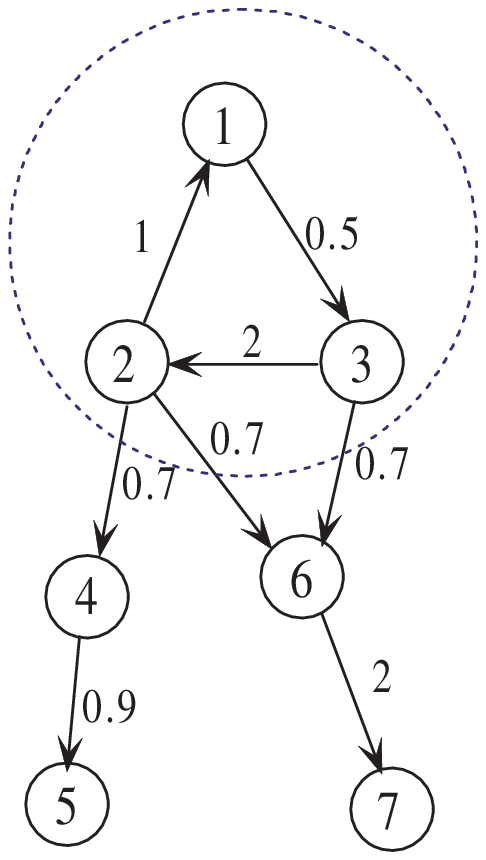}
      \includegraphics[scale=0.8]{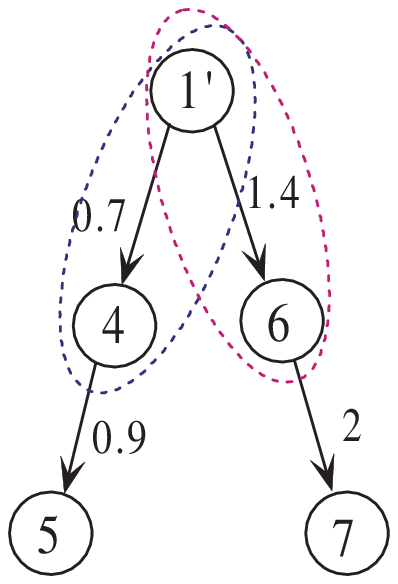}
      \includegraphics[scale=0.8]{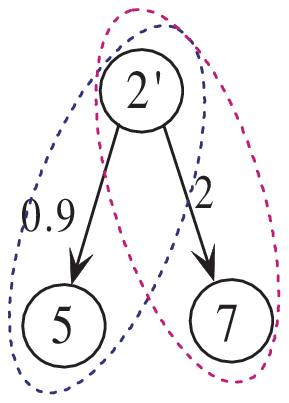}

\ \ \ \ \  \ \ \ \  \ \      (a)\ \ \ \ \ \ \ \  \ \ \ \ \ \  \ \ \
\ \ \ \  \ \ \ (b)\ \ \ \  \ \ \  \ \ \ \ \ \ \  \ \ \ \ \ \ \ (c)
\small

Suppose that (a) is the underlying communication topology, where the
numbers in the circles are the indices of agents, and the numbers
near the edges are the weights. The vertices in the dashed circle in
(a) are the leaders, the local communication topology among which is
strongly connected. After a period of time, the states of them will
agree, and they can be viewed as one agent. The equivalent topology
is (b), where $1'$ represents the virtue agent. In (b), by the
conclusion of the second step, agent $4$ and agent $6$ will also
agree with the leaders separately after a time. And then the
equivalent topology becomes (c). And finally, they will reach a
consensus on states.

 \caption{Demonstration of the third step of the proof of Theorem
\ref{dl1}}
      \label{demostep}
\end{figure}

{\it Step 3:} Finally, suppose that the communication topology
$\mathcal{G}(\A)$ has a spanning tree.

Consider  another directed graph $\G$, whose vertex set is the set
of all strongly connected components of $\G(\A)$, denoted by
$\{u_1,u_2,\dots,u_m\}$, $m\leq n$. $(u_i,u_j)\in\mathcal{E}(\G)$ if
and only if there exist $v_{i'}\in\mathcal{V}(u_i)$ and
$v_{j'}\in\mathcal{V}(u_j)$ such that $(v_{i'},
v_{j'})\in\mathcal{E}(\G(\A))$. Then $\G$ is uniquely determined by
$\A$ and is a directed tree. Denote the root of $\G$ by $u_{s_1}$,
where its vertex set corresponds to the leader set. Consider the
states of agents associated with $\mathcal{E}(u_i)$. Since $\G$ is a
directed tree, there exists a path $u_{s_1}, u_{s_2}, \dots,
u_{s_l}=u_i$ from the root $u_{s_1}$ to $u_i$. For the agents
corresponding to the vertex set of $u_{s_1}$, since the local
communication topology $u_{s_1}$ among them is strongly connected,
and the dynamics of them is not affected by others, by the
conclusion of the first step, their states will reach a consensus in
finite time $t_1$. After time $t_1$, the states of those agents are
time-invariant. Therefore, under protocol \eqref{pro1}, they can be
viewed as one agent for the other agents. Denote this virtue agent
by $\bar{v}_1$. Define a new  weighted directed graph $\G_1$. The
vertex set is $\{\bar{v}_1\}\cup\mathcal{E}(u_{s_2})$. For any
$v_j\in\mathcal{E}(u_{s_2})$, $(\bar{v}_1, v_j)\in
\mathcal{E}(\G_1)$ if and only if there exists
$v_k\in\mathcal{E}(u_{s_1})$ such that
$(v_k,v_j)\in\mathcal{E}(\G(\A))$. And the weight of $(\bar{v}_1,
v_j)$ is the sum of the weights of all such edges as $(v_k,v_j)$.
The induced weighted subgraph of $\G_1$ by the vertices of $u_{s_2}$
is also $u_{s_2}$ and inherits its weights in $\G(\A)$. Notice that
the dynamics of agents associated with $\mathcal{E}(u_{s_2})$ under
topology $\G(\A)$ is the same as the dynamics under topology $\G_1$.
Hence, by the conclusion of the second step, there exists time
$t_2\geq t_1$, at which those agents will reach a consensus on
states with the leaders. By induction, the states of agents
corresponding to $\mathcal{E}(u_j)$ will be  the same as the
leaders' final states in a finite time. The illustration of the
proof is shown in Fig. \ref{demostep}.

To conclude, the system solves a finite-time consensus problem.

\end{proof}

From the first step of the proof of Theorem \ref{dl1}, we can see
that
\begin{tl}
Suppose that $\G(\A)$ is strongly connected and there exists
$\w\in\R^n$ such that $\w>0$ and $\w^TL(\A)=0$. If the following
protocol
\begin{equation}\label{pro3}
    u_i=\sum_{j=1}^na_{ij}(\sgn{x_j}{\alpha_j}-\sgn{x_i}{\alpha_i}),i\in\n,
\end{equation}
where $0<\alpha_i,\alpha_j<1$, is taken and the initial state
$\x(0)$ satisfies $\w^T\x(0)=0$, then the states of agents will
reach zero in finite time.
\end{tl}

{\it Remark:} Notice that
$\w^T\dot{\x}(t)=\w^TL(\A)\sgn{\x}{\al}=0$. Then $\w^T\x(t)$ is
time-invariant and $\w^T\x(t)\equiv 0$. Thus the above system is the
same as the system of $\y$ in the first step of the proof of Theorem
\ref{dl1}. And it is easy to check that the above system satisfies
the Lipschitz condition at all points in $\spans(\yi)$ except the
origin, and thus it is impossible for the states of agents reach
consensus in finite time on other agreement points.

\subsection{Networks Under Protocol \eqref{pro2}}

In this subsection,  we present conditions under which protocol
\eqref{pro2} solves a finite-time consensus problem. In fact when
the topology is dynamically changing, protocol \eqref{pro2} is also
applicable. This case will be discussed thoroughly in the subsequent
subsection.

Unlike the discussions on protocol \eqref{pro1}, here will present
several theorems, which are valid in different cases, and are of
interest themselves.
\begin{dl}\label{dl2}
Suppose that communication topology $\mathcal{G}(\A)$ is undirected
and connected. Then protocol \eqref{pro2} solves the finite-time
average-consensus  problem, if $\alpha_{ij}=\alpha_{ji}$ for all
$i\in\n$, $j\in\mathcal{N}(\G(\A),i)$.
\end{dl}

\begin{proof}
First, complete definition of $\alpha_{ij}$ is provided. For any
$i\in\n$, if $j\not\in\N(\G(\A),i)$, set
$\alpha_{ij}=r\leq\max_{k\in\n,l\in\mathcal{N}(\G(\A),k)}\alpha_{kl}$.
Clearly in this case $a_{ij}=0$.

 Since $a_{ij}=a_{ji}$ and $\alpha_{ij}=\alpha_{ji}$, for all $i,j\in \n$,
we get that
\begin{equation}\label{eqdl2eq1}
    \sum_{i=1}^n\dot{x}_i(t)=0.
\end{equation}
Let
\[\kappa=\frac{1}{n}\sum_{i=1}^n x_i(t).\] By \eqref{eqdl2eq1}, $\kappa$
is time-invariant. Let $x_i(t)=\kappa+\delta_i(t)$ and let
$\de(t)=[\delta_1(t),\delta_2(t),\dots,\delta_n(t)]^T$. Then
$\dot{\delta}_i(t)=\dot{x}_i(t)$. In \cite{R. Olfati-Saber and R. M.
Murray 1}, $\de(t)$ is referred to as the {\it group disagreement
vector}. Consequently,
$\yi^T\de(t)=\sum_{i=1}^nx_i(t)-n\kappa=\sum_{i=1}^nx_i(t)-\sum_{i=1}^nx_i(t)=0$.
We take Lyapunov function
\[V_3(\de(t))=\frac{1}{2}\sum_{i=1}^n \delta_i^2(t).\]
 The remainder of
the proof is to show that $V_3(t)$ satisfies the conditions of Lemma
\ref{ylfinitetimeV} for some $K$ and $\alpha$ and then $V_3(t)$
reaches zero in finite time, which yields that all agents' states
 reach a consensus in finite time and the final state is $\kappa$,
the average of their initial states.

Differentiate $V_3(t)$ with respect to $t$.
\begin{align}
\frac{dV_3(t)}{dt}=&\sum_{i=1}^n \delta_i(t) \dot{\delta}_i(t)\nonumber\\
=&\frac{1}{2}\sum_{i,j=1}^n\left(a_{ij}\delta_i
\sgn{\delta_j-\delta_i}{\alpha_{ij}}  +a_{ji}\delta_j
\sgn{\delta_i-\delta_j}{\alpha_{ji}}\right)\label{eq1dl2}\\
=&\frac{1}{2}\sum_{i,j=1}^n a_{ij}(\delta_i-\delta_j)
\sgn{\delta_j-\delta_i}{\alpha_{ij}}\nonumber\\
=&-\frac{1}{2}\sum_{i,j=1}^na_{ij}|\delta_j-\delta_i|^{1+\alpha_{ij}}\nonumber\\
=&-\frac{1}{2}\sum_{i,j=1}^n\left(a_{ij}^{\frac{2}{1+\alpha_0}}\left((\delta_j-\delta_i)^{2}\right)^{\frac{1+\alpha_{ij}}{1+\alpha_0}}\right)^{\frac{1+\alpha_0}{2}},\nonumber
\end{align}
where $\alpha_0=\max_{ij}\alpha_{ij}$ and equation \eqref{eq1dl2}
follows from that $\delta_j-\delta_i=x_j-x_i$. From $0<\alpha_0<1$,
$\frac{1}{2}<\frac{1+\alpha_0}{2}<1$.

Suppose that $V_3(t)\not=0$. By Lemma \ref{yl1},
\begin{align}
\frac{dV_3(t)}{dt}& \leq -\frac{1}{2} \left(\sum _{i,j=1}^n
a_{ij}^{\frac{2}{1+\alpha_0}}((\delta_i-\delta_j)^{2})^{\frac{1+\alpha_{ij}}{1+\alpha_0}}\right)^\frac{1+\alpha_0}{2}\nonumber\\
&=  -\frac{1}{2}\left(\frac{\sum _{i,j=1}^n
a_{ij}^{\frac{2}{1+\alpha_0}}((\delta_i-\delta_j)^{2})^{\frac{1+\alpha_{ij}}{1+\alpha_0}}}{\sum
_{i,j=1}^na_{ij}^{\frac{2}{1+\alpha_0}}(\delta_i-\delta_j)^{2}}\frac{\sum
_{i,j=1}^na_{ij}^{\frac{2}{1+\alpha_0}}(\delta_i-\delta_j)^{2}}{V_3(t)}V_3(t)\right)^\frac{1+\alpha_0}{2}\label{eqdl2eq2}
\end{align}

The last equation follows from that $\sum
_{i,j=1}^na_{ij}^{\frac{2}{1+\alpha_0}}(\delta_i-\delta_j)^{2}\not=0$.
In fact, if $\sum
_{i,j=1}^na_{ij}^{\frac{2}{1+\alpha_0}}(\delta_i-\delta_j)^{2}=0$,
then by the connectivity of $\G(\A)$, $\delta_i=\delta_j$ for all
$i,j\in\n$, namely, $\de\in\spans(\yi)$.  Because $\yi^T\de=0$,
$\de=0$, and thus $V_3(t)=0$, which contradicts our assumption.
Therefore, $\sum
_{i,j=1}^na_{ij}^{\frac{2}{1+\alpha_0}}(\delta_i-\delta_j)^{2}\not=0$.

Next, we give the lower bounds of the first two qualities in the big
brackets of  equation \eqref{eqdl2eq2}.

By Property \ref{xz1}, $\max_i x_i-\min_i x_i$ is non-increasing.
For any $i,j\in\n$, $|\delta_i(t)-\delta_j(t)|\leq \max_k
x_k(t)-\min_k x_k(t)\leq \max_k x_k(0)-\min_k x_k(0)$. Let
\[K_4=\frac{1}{\sum_{i,j=1}^n a_{ij}^{\frac{2}{1+\alpha_0}}}
 \min_{i,j\in\n\atop a_{ij}\not=0}a_{ij}^{\frac{2}{1+\alpha_0}}(\max_kx_k(0)-\min_k
x_k(0))^{2(\frac{1+\alpha_{ij}}{1+\alpha_0}-1)},\] which is
positive. Let $(i_0,j_0)=\arg\max_{i,j\in\n\atop
a_{ij}\not=0}(\delta_i-\delta_j)^{2}$. Since
$\frac{1+\alpha_{ij}}{1+\alpha_0}-1\leq 0$, with the  same arguments
as in the proof of inequality \eqref{eqdl1inq2},
\begin{align*}
\frac{\sum _{i,j=1}^n
a_{ij}^{\frac{2}{1+\alpha_0}}((\delta_i-\delta_j)^{2})^{\frac{1+\alpha_{ij}}{1+\alpha_0}}}{\sum
_{i,j=1}^na_{ij}^{\frac{2}{1+\alpha_0}}(\delta_i-\delta_j)^{2}}&\geq
\frac{
a_{i_0j_0}^{\frac{2}{1+\alpha_0}}((\delta_{i_0}-\delta_{j_0})^{2})^{\frac{1+\alpha_{i_0j_0}}{1+\alpha_0}}}{\left(\sum_{i,j=1}^n
a_{ij}^{\frac{2}{1+\alpha_0}}\right)(\delta_{i_0}-\delta_{j_0})^{2}}\\
&\geq K_4.
\end{align*}

Now consider the second quality. Let $\B=[b_{ij}]\in\R^{n\times n}$,
where $b_{ij}=a_{ij}^{\frac{2}{1+\alpha_0}}$. Then by Lemma
\ref{yl2},
\[\sum _{i,j=1}^n
a_{ij}^{\frac{2}{1+\alpha}}(\delta_i-\delta_j)^2=\sum_{i,j=1}^nb_{ij}(\delta_i-\delta_j)^2=2\de^T
L(\B)\de,\] and since $\de\perp\yi$,
\[\frac{\sum
_{i,j=1}^na_{ij}^{\frac{2}{1+\alpha_0}}(\delta_i-\delta_j)^{2}}{V_3(t)}=\frac{2\de^TL(\B)\de}{\frac{1}{2}\de^T\de}\geq
4\lambda_2(L(\B))>0.\] Therefore,
\begin{align*}
    \frac{d V_3(t)}{dt}\leq
    -\frac{1}{2}\left(4
    K_4\lambda_2(L(\B))\right)^{\frac{1+\alpha_0}{2}}V_3(t)^{\frac{1+\alpha_0}{2}}.
\end{align*}

Thus, by Lemma \ref{ylfinitetimeV}, system \eqref{xt} solves the
finite-time average-consensus problem.

\end{proof}

The next theorem is on the case with a leader.

\begin{dl}\label{dlpro2leader}
Suppose that the system consists of one leader and $n-1$ followers
and the local communication topology among the followers is
undirected and connected. If $\alpha_{ij}=\alpha_{ji}$ for all $i,j$
such that agent $i$ and agent $j$ are neighbors of each other, then
protocol \eqref{pro2} solves a finite-time consensus problem and the
final common state is the state of the leader.
\end{dl}
\begin{proof}
First we define the undefined $\alpha_{ij}$ as in the proof of
Theorem \ref{dl2} and let $m=n-1$. Without loss of generality,
suppose that agents $1,2,\dots,m$ are the followers and agent $n$ is
the leader. From the provided conditions of this theorem, we have
that $a_{n1}=a_{n2}=\dots=a_{nn}=0$. Let $\bar{\A}=[a_{ij}]_{1\leq
i,j\leq m}\in\R^{m \times m}$, let
$\bar{\bb}=[b_1,b_2,\dots,b_m]^T=[a_{1n},a_{2n},\dots, a_{mn}]^T$
and let $\alpha_0=\max_{ij} \alpha_{ij}$. Then $\bar{\A}^T=\bar{\A}$
and $\bar\bb\not=0$. Rewrite protocol \eqref{pro2}
\[u_i=\sum_{j=1}^ma_{ij}\sgn{x_j-x_i}{\alpha_{ij}}+ b_i\sgn{x_n-x_i}{\alpha_{in}}, i\in\m.\]
Let $\delta_i=x_i-x_n$, $i\in\n$. Then $\delta_n\equiv0$ and
\[\dot{\delta_i}=\dot{x}_i=\sum_{j=1}^ma_{ij}\sgn{\delta_j-\delta_i}{\alpha_{ij}}- b_i\sgn{\delta_i}{\alpha_{in}}, i\in\m.\]

Also consider the Lyapunov function
$V_3(t)=\frac{1}{2}\sum_{i=1}^n\delta_i^2$.
\begin{align*}
  \frac{d
  V_3(t)}{dt}=&\sum_{i=1}^m\delta_i\left(\sum_{j=1}^ma_{ij}\sgn{\delta_j-\delta_i}{\alpha_{ij}}-
  b_i\sgn{\delta_i}{\alpha_{in}}\right)\\
  =&-\frac{1}{2}\sum_{i,j=1}^ma_{ij}|\delta_j-\delta_i|^{1+\alpha_{ij}}-\sum_{i=1}^mb_i|\delta_i|^{1+\alpha_{in}}\\
  =&-\frac{1}{2}\sum_{i,j=1}^m
  \left(a_{ij}^{\frac{2}{1+\alpha_0}}|\delta_j-\delta_i|^{\frac{2(1+\alpha_{ij})}{1+\alpha_0}}\right)^{\frac{1+\alpha_0}{2}}-
   \sum_{i=1}^{m}\left(b_i^{\frac{2}{1+\alpha_0}}|\delta_i|^{\frac{2(1+\alpha_{in})}{1+\alpha_0}}\right)^{\frac{1+\alpha_0}{2}}\\
   \leq &-\frac{1}{2}\left( \sum_{i,j=1}^ma_{ij}^{\frac{2}{1+\alpha_0}}|\delta_j-\delta_i|^{\frac{2(1+\alpha_{ij})}
   {1+\alpha_0}}+ 2\sum_{i=1}^{m}b_i^{\frac{2}{1+\alpha_0}}|\delta_i|^{\frac{2(1+\alpha_{in})}{1+\alpha_0}}
    \right)^{\frac{1+\alpha_0}{2}}.
\end{align*}
The last inequality is obtained by Lemma \ref{yl1}.

For simplicity, let
$G_1(\de)=\sum_{i,j=1}^ma_{ij}^{\frac{2}{1+\alpha_0}}|\delta_j-\delta_i|^{\frac{2(1+\alpha_{ij})}
   {1+\alpha_0}}+
   2\sum_{i=1}^{m}b_i^{\frac{2}{1+\alpha_0}}|\delta_i|^{\frac{2(1+\alpha_{in})}{1+\alpha_0}}$
   and let $G_2(\de)=\sum_{i,j=1}^ma_{ij}^{\frac{2}{1+\alpha_0}}(\delta_j-\delta_i)^{2}+
   2\sum_{i=1}^{m}b_i^{\frac{2}{1+\alpha_0}}\delta_i^{2}$.
With
   the same arguments as in the proof of inequality
   \eqref{eqdl1inq2},
   \[\frac{G_1(\de)}{G_2(\de)}\geq K_4,\]
   where $K_4$ is defined in Theorem \ref{dl2}.

Let $\bar{\B}=[a_{ij}^{\frac{2}{1+\alpha_0}}]_{1\leq i,j\leq m}$ and
let
$\tilde\bb=[b_1^{\frac{2}{1+\alpha_0}},b_2^{\frac{2}{1+\alpha_0}},\dots,b_m^{\frac{2}{1+\alpha_0}}]^T$.
By Corollary \ref{tlLb}, $L(\bar{\B})+\diag(\tilde\bb)$ is positive
definite. Its smallest eigenvalue is denoted by
$\lambda_1(L(\bar{\B})+\diag(\tilde\bb))$. Then
\[\frac{G_2(\de)}{V_3(t)}=\frac{2[\delta_1,\delta_2,\dots,\delta_m]^T(L(\bar{\B})+\diag(\tilde\bb))[\delta_1,\delta_2,\dots,\delta_m]}
{\frac{1}{2}[\delta_1,\delta_2,\dots,\delta_m]^T[\delta_1,\delta_2,\dots,\delta_m]}\geq
4\lambda_1(L(\bar{\B})+\diag(\tilde\bb)).\]

Let $\B$ be the same as in the proof of Theorem \ref{dl2}. We have
\[\B=\left[
       \begin{array}{cc}
         \bar\B & \tilde\bb \\
         0 & 0 \\
       \end{array}
     \right].
\]
Therefore,
\[L(\B)=\left[
       \begin{array}{cc}
         L(\bar\B)+ \diag(\tilde\bb)& -\tilde\bb \\
         0 & 0 \\
       \end{array}
     \right].\]
Its eigenvalues are all nonnegative numbers, and the second smallest
eigenvalue of $L(\B)$, denoted by $\lambda_2(L(\B))$, is
$\lambda_1(L(\bar{\B})+\diag(\tilde\bb))$. Thus,
\[\frac{d V_3(t)}{dt}\leq -\frac{1}{2}\left(4K_4\lambda_2(L(\B))\right)^{\frac{1+\alpha_0}{2}}V_3(t)^{\frac{1+\alpha_0}{2}}.\]
By Lemma \ref{ylfinitetimeV}, $V_3(t)$ will reach zero in finite
time, which implies that all agents' states reach the common value
$x_n(0)$ in a finite time.
\end{proof}

In the above analysis, we only consider the case when $\G(\A)$ is
undirected or the local communication topology among the followers
is undirected. To some extent, the above two results can be extended
to the case when the communication topology belongs to the set of
some special kinds of directed graphs.

\begin{dy}[\cite{chu detailbalance}]
Communication topology $\G(\A)$ is said to satisfy the {\it detail
balance condition} in weights if there exist some scalars
$\omega_i>0$, $i\in\n$, such that $\omega_ia_{ij}=\omega_ja_{ji}$
for all $i,j\in\n$.
\end{dy}

We can see that in the detail-balanced communication topology, the
communication channels among agents are also bidirectional like in
undirected graphs, but with different weights.

\begin{tl}\label{tldl2}
Protocol \eqref{pro2} solves a finite-time consensus problem if
$\G(\A)$ is strongly connected and  detail-balanced, and
$\alpha_{ij}=\alpha_{ji}$ for all $i\in\n$, $j\in\bN(\G(\A),i)$.
\end{tl}
\begin{proof}
Let $\w=[\omega_1,\omega_2,\dots,\omega_n]^T$ such that
$\omega_ia_{ij}=\omega_ja_{ji}$ for all $i,j\in\n$ and define the
undefined $\alpha_{ij}$ as in the proof of Theorem \ref{dl2}. It is
easy to check that $\diag(\w)\A$ is symmetric and
\[\sum_{i=1}^n\omega_i\dot{x}_i(t)=0.\]

Let $\kappa=\frac{1}{\sum_{i=1}^n
\omega_i}\sum_{i=1}^n\omega_ix_i(t)$.  Then $\kappa$ is
time-invariant. Let $\delta_i(t)=x_i(t)-\kappa$. Then
$\w^T\de(t)=0$. We take Lyapunov function
$V_4(t)=\frac{1}{2}\sum_{i=1}^n\omega_i\delta_i(t)^2$.

\[\frac{dV_4(t)}{dt}=-\frac{1}{2}\sum_{i,j=1}^n\omega_ia_{ij}|\delta_j-\delta_i|^{1+\alpha_{ij}}.\]
With the same argument as in the proof of Theorem \ref{dl2}, there
exists $K_5>0$ such that
\[\frac{dV_4(t)}{dt}\leq -K_5 V_4(t)^{\frac{1+\alpha_0}{2}},\]
where $\alpha_0=\max_{ij}\alpha_{ij}$. The details are omitted. And
thus this system solves a finite-time consensus problem, and the
final state is
$\frac{1}{\sum_{i=1}^n\omega_i}\sum_{i=1}^n\omega_ix_i(0)$, which
can be viewed as a {\it weighted average-consensus function}.
\end{proof}

The following corollary is also obtainable but the proof is omitted.
\begin{tl}\label{tldlpro2leader}
Suppose that the conditions are the same as in Theorem
\ref{dlpro2leader}, except that the local communication topology
among the followers is strongly connected and detail-balanced. Then
protocol \eqref{pro2} solves a finite-time consensus problem.
\end{tl}

Based on the results of Corollary \ref{tldl2} and Corollary
\ref{tldlpro2leader}, we can still extend those results further,
though, for protocol \eqref{pro2}, we cannot obtain results as
desirable as for protocol \eqref{pro1}.
\begin{dl}
If $\G(\A)$ has a spanning tree and each strongly connected
component is detail-balanced, and $\alpha_{ij}=\alpha_{ji}$ for all
$i,j$ such that agent $i$ and agent $j$ are neighbors of each other,
then protocol \eqref{pro2} solves a finite-time consensus problem.
\end{dl}

The proof of the above theorem is the same as the third step of the
proof of Theorem \ref{dl1}.

\subsection{Performance Analysis}
In this subsection, we investigate the relationship between the
convergence time and the other factors, such as the underlying
communication topology, parameters $\alpha_i$ and $\alpha_{ij}$ in
the proposed protocols, and the initial states. The {\it convergence
time } is defined as the time that the system spends to reach a
consensus. For the studied nonlinear system, the precise time of
reaching the consensus state is hard to given. However, we can
analyze it by investigating its upper bound obtained by Lemma
\ref{ylfinitetimeV}.

We claim that $V_1(t), V_2(t), V_3(t)$ and $V_4(t)$  measure how
much agents' states differ from each other. $ V_3(t)$ and $V_4(t)$
clearly measure the disagreement of agents' current states with
their common final state by their definitions. Now consider
$V_1(t)$. We take the same assumptions and notations as in Theorem
\ref{dl1}. Let $\w^{\bot}=\{\bxi\in\R^n:\w^T\bxi=0\}$. Then
$\w^{\bot}$ is $L(\A)$-invariant. Furthermore, for any
$\bxi\in\w^{\bot}$, $L(\A)\bxi=0\iff \bxi=0$. Therefore,
$\|\bxi\|_{L(\A)}=\|L(\A)\bxi\|_2$ defines a vector norm on
$\w^{\bot}$ (see \cite{Horn}, Theorem 5.3.2). Furthermore, for any
$\x$, let $\x=r\yi+\x_2$, where $r\in\R$ and $\x_2\in\w^{\perp}$.
\[\|L(\A)\x\|_{2}=\sqrt{\x^TL(\A)^TL(\A)\x}=\sqrt{\x_2^TL(\A)^TL(\A)\x_2}=\|\x_2\|_{L(\A)}.\]
 And apparently, $V_1(t)$ is a
measure of the length of $L(\A)\x$. Especially, if all $\alpha_i$
are equal to $\alpha$, $0<\alpha<1$, and $\w=\yi$,
$V_1(\x(t))=\frac{1}{1+\alpha}(\|L(\A)\x\|_{(1+\alpha)})^{1+\alpha}$.
Therefore, $V_1(t)$ can be seen as a measure of the length of vector
$\x_2$, which reflects the disagreement of $\x$ with subspace
$\spans(\yi)$. In the same way, we can show that $V_2(t)$ also
measures the disagreement of agents' states. Undoubtedly, larger
initial value of $V_1(t), V_2(t), V_3(t)$ and $V_4(t)$ will result
in longer convergence time of the considered system by Lemma
\ref{ylfinitetimeV}. This fact also can be reflected in the
estimation of $K_2$, $K_3$ and $K_4$.

It is well known that, if $\alpha_i$ or $\alpha_{ij}=1$ for all
$i,j$, then the system at most asymptotically solves a consensus
problem, and thus the convergence time becomes infinitively long. In
our model, since $\lim_{\alpha_i\to
1}\frac{2\alpha_0}{1+\alpha_0}=1$ and $\lim_{\alpha_{ij}\to
1}\frac{1+\alpha_0}{2}=1$, by Lemma \ref{ylfinitetimeV}, the
estimated upper bound of convergence time
 tends to infinity, if all $\alpha_i$ or $\alpha_{ij}$ tend to $1$.

In \cite{R. Olfati-Saber and R. M.  Murray 1}, it was proved that
the larger the algebraic connectivity of the underlying
communication topology is, the faster the system converges. To
illustrate the dependency  of the convergence time under protocol
\eqref{pro1} or \eqref{pro2} on the algebraic connectivity of the
communication topology, we assume that $\G(\A)$ is undirected and
connected, and all $\alpha_i$ or $\alpha_{ij}$ are all equal to
$\alpha$, $0<\alpha<1$. For protocol \eqref{pro1}, we consider
Lyapunov function
$V_5(t)=\frac{1}{4}\sum_{i,j=1}^na_{ij}(x_j(t)-x_i(t))^2$ instead of
$V_1(t)$. It can be shown that $V_5(t)=0\iff \x\in\spans(1)$ and
\begin{equation}\label{ineqv5-a}
    \frac{dV_5(t)}{dt}\leq - (2\lambda_2(L(\A)))^{\frac{1+\alpha}{2}}V_5(t)^{\frac{1+\alpha}{2}}.
\end{equation}
For protocol \eqref{pro2}, $K_4$ can be taken to be $1$, and
\begin{equation}\label{ineqv3-a}
    \frac{dV_3(t)}{dt}\leq-\frac{1}{2}(4\lambda_2(L(\B)))^{\frac{1+\alpha}{2}}V_3(t)^{\frac{1+\alpha}{2}},
\end{equation}
where $\B=[a_{ij}^{\frac{2}{1+\alpha}}].$

The above two differential inequalities show that larger algebraic
connectivity of $\G(\A)$ or $\G(\B)$ can lead to shorter the
convergence time (see Lemma \ref{ylfinitetimeV}).

The following two theorems compare the convergence rates of protocol
\eqref{pro1} or protocol \eqref{pro2} with different parameters
$\alpha_i$ or $\alpha_{ij}$ under the same communication topology
respectively. It is shown that smaller $\alpha_i$ or $\alpha_{ij}$
can get better convergence rate when agents states differ a little
from each other, and larger $\alpha_i$ or $\alpha_{ij}$ can get
better convergence rate when agents states differ a lot from each
other. Therefore, in order to get better convergence rate, we can
change the values of protocol parameters $\alpha_i$ or $\alpha_{ij}$
based on the states of agents as system evolves.

\begin{dl}\label{dlcomparea-1}
Suppose that $\G(\A)$ is undirected and connected, protocol
\eqref{pro1} is applied, and all $\alpha_i$ are  equal, denoted by
$\alpha$. Then the eigenvalues of $L(\A)$ are nonnegative numbers
(by Lemma \ref{yl2}) and let them be $\lambda_1(L(\A)),
\lambda_2(L(\A))$, $\dots$, $\lambda_n(L(\A))$ in the increasing
order. Since $\G(\A)$ is connected, $\lambda_1(L(\A))=0$ and
$\lambda_2(L(\A))>0$.
 Given
$0<\alpha_*<\alpha^*<1$, let
${\epsilon}^*=\frac{1}{2\lambda_2(L(\A))}n^{\frac{1-\alpha_*}{\alpha^*-\alpha_*}}$
and let ${\epsilon}_*= \frac{1}{2\lambda_n(L(\A))}$. Consider
Lyapunov function
$V_5=\frac{1}{4}\sum_{i,j=1}^na_{ij}(x_j(t)-x_i(t))^2$. If
$V_5(t)\geq {\epsilon}^*$, then
\[\left.\frac{dV_5(t)}{dt}\right|_{\alpha=\alpha^*}\leq\left.\frac{dV_5(t)}{dt}\right|_{\alpha=\alpha_*},\]
and if $V_5(t)\leq {\epsilon}_*$,
\[\left.\frac{dV_5(t)}{dt}\right|_{\alpha=\alpha_*}\leq\left.\frac{dV_5(t)}{dt}\right|_{\alpha=\alpha^*}.\]
\end{dl}
\begin{proof}
Since
\begin{equation}\label{eqVdifferentapro1-1}
    \frac{dV_5(t)}{dt}=-\sum_{i=1}^n\left|\sum_{j=1}^na_{ij}(x_j-x_i)\right|^{1+\alpha}
=-\sum_{i=1}^n\left(\left(\sum_{j=1}^na_{ij}(x_j-x_i)\right)^2\right)^{\frac{1+\alpha}{2}},
\end{equation}
by Lemma \ref{yl1},
\begin{equation}\label{eqVdifferentapro1-2}
    -n^{\frac{1-\alpha}{2}}\left(\sum_{i=1}^n\left(\sum_{j=1}^na_{ij}(x_j-x_i)\right)^2\right)^{\frac{1+\alpha}{2}}\leq\frac{dV_5(t)}{dt}
\leq-\left(\sum_{i=1}^n\left(\sum_{j=1}^na_{ij}(x_j-x_i)\right)^2\right)^{\frac{1+\alpha}{2}}.
\end{equation}

It is easy to see that
\[\sum_{i=1}^n\left(\sum_{j=1}^na_{ij}(x_j-x_i)\right)^2=\x^TL(\A)^TL(\A)\x,\]
and \[V_5(t)=\frac{1}{2}\x^TL(\A)\x.\]
 Let \[\D=\left[
         \begin{array}{cccc}
           0 &  &  &  \\
            & \lambda_2(L(\A)) &  &  \\
            &  &  \ddots&  \\
            &  &  &  \lambda_n(L(\A))\\
         \end{array}
       \right].\] Since  $L(\A)$ is symmetric,
there exists an orthogonal matrix $\T\in\R^{n\times n}$ such that
\[L(\A)=\T^T\D\T.
\]
And,
\begin{align*}
\frac{\x^T
L(\A)^TL(\A)\x}{\frac{1}{2}\x^TL(\A)\x}=\frac{2\x^T\T^T\D\T\T^T\D\T\x}{\x^T\T^T\D\T\x}=\frac{2\x^T\T^T\D^2\T\x}{\x^T\T^T\D\T\x}.
\end{align*}
Therefore,
\begin{equation}\label{eqVdifferentapro1-3}
    2\lambda_2(L(\A))\leq \frac{\x^T
L(\A)^TL(\A)\x}{\frac{1}{2}\x^TL(\A)\x}\leq 2\lambda_n(L(\A)).
\end{equation}

If $V_5(t)\geq{\epsilon}^*$, then
\begin{align*}
    \left.\frac{dV_5(t)}{dt}\right|_{\alpha=\alpha^*}\leq &-\left(\x^TL(\A)^TL(\A)\x\right)^{\frac{1+\alpha^*}{2}}\mbox{( by \eqref{eqVdifferentapro1-2})}\\
    =-&
    \left(\x^TL(\A)^TL(\A)\x\right)^{\frac{1+\alpha_*}{2}}\left(\x^TL(\A)^TL(\A)\x\right)^{\frac{\alpha^*-\alpha_*}{2}}\\
    \leq &-\left(\x^TL(\A)^TL(\A)\x\right)^{\frac{1+\alpha_*}{2}}\left(2\lambda_2(L(\A))V_5(t)
    \right)^{\frac{\alpha^*-\alpha_*}{2}}\mbox{( by \eqref{eqVdifferentapro1-3})}\\
    \leq
    &-n^{\frac{1-\alpha_*}{2}}\left(\x^TL(\A)^TL(\A)\x\right)^{\frac{1+\alpha_*}{2}} \mbox{ (by $V_5(t)\geq \epsilon^*$)}\\
    \leq&\left.\frac{dV_5(t)}{dt}\right|_{\alpha=\alpha_*}\mbox{( by \eqref{eqVdifferentapro1-2})}.
\end{align*}

If $V_5(t)\leq{\epsilon}_*$, by \eqref{eqVdifferentapro1-3},
\[\sum_{i=1}^n\left(\sum_{j=1}^na_{ij}(x_j-x_i)\right)^2\leq 2\lambda_n(L(\A))V_5(t)\leq 1.\]
Therefore, for any $i\in\n$, $
\left(\sum_{j=1}^na_{ij}(x_j-x_i)\right)^2\leq 1$. And because
$\alpha_*<\alpha^*$, $\frac{1+\alpha_*}{2}< \frac{1+\alpha^*}{2}$.
Thus
\[\left(\left(\sum_{j=1}^na_{ij}(x_j-x_i)\right)^2\right)^{\frac{1+\alpha_*}{2}}
\geq
\left(\left(\sum_{j=1}^na_{ij}(x_j-x_i)\right)^2\right)^{\frac{1+\alpha^*}{2}}.\]
By \eqref{eqVdifferentapro1-1},
\[\left.\frac{dV_5(t)}{dt}\right|_{\alpha=\alpha_*}\leq\left.\frac{dV_5(t)}{dt}\right|_{\alpha=\alpha^*}.\]
\end{proof}

\begin{dl}\label{dlcomparea-2}
Suppose that $\G(\A)$ is undirected and connected, protocol
\eqref{pro2} is applied, and all $\alpha_{ij}$ are equal, denoted by
$\alpha$. Let $\B=[a_{ij}^{\frac{2}{1+\alpha}}]$.  Then the
eigenvalues of $L(\B)$ are nonnegative numbers (by Lemma \ref{yl2})
and let the second smallest and the largest eigenvalues be
$\lambda_2(L(\B)), \lambda_n(L(\B))$ respectively (which are
positive).
 Given
$0<\alpha_*<\alpha^*<1$, let
\[\epsilon^*=2^{-2}n^{\frac{2(1-\alpha_*)}{\alpha^*-\alpha_*}}
\left(\lambda_n(L(\B|_{\alpha=\alpha_*}))\right)^{\frac{1+\alpha_*}{\alpha^*-\alpha_*}}
\left(\lambda_2(L(\B|_{\alpha=\alpha^*}))\right)^{\frac{1+\alpha^*}{\alpha_*-\alpha^*}},\]
and
 let
\[\epsilon_*= 2^{-2}n^{\frac{2(1-\alpha^*)}{\alpha_*-\alpha^*}}
\left(\lambda_n(L(\B|_{\alpha=\alpha^*}))\right)^{\frac{1+\alpha^*}{\alpha_*-\alpha^*}}
\left(\lambda_2(L(\B|_{\alpha=\alpha_*}))\right)^{\frac{1+\alpha_*}{\alpha^*-\alpha_*}}.\]
Consider Lyapunov function
$V_3(\de(t))=\frac{1}{2}\sum_{i=1}^n\delta_i^2(t)$ (where $\de(t)$
was defined in the proof of Theorem \ref{dl2}). If $V_3(t)\geq
\epsilon^*$, then
\[\left.\frac{dV_3(t)}{dt}\right|_{\alpha=\alpha^*}\leq\left.\frac{dV_3(t)}{dt}\right|_{\alpha=\alpha_*},\]
and if $V_3(t)\leq {\epsilon}_*$,
\[\left.\frac{dV_3(t)}{dt}\right|_{\alpha=\alpha_*}\leq\left.\frac{dV_3(t)}{dt}\right|_{\alpha=\alpha^*}.\]
\end{dl}
\begin{proof}
In the proof of Theorem \ref{dl2}, we get that
\[\frac{dV_3(t)}{dt}=-\frac{1}{2}\sum_{i,j=1}^n\left(a_{ij}^{\frac{2}{1+\alpha}}
(\delta_j-\delta_i)^{2}\right)^{\frac{1+\alpha}{2}},\] and by Lemma
\ref{yl1},
\begin{equation}\label{eqdlcomparea-2-1}
    -\frac{n^{1-\alpha}}{2}\left(\sum_{i,j=1}^na_{ij}^{\frac{2}{1+\alpha}}
(\delta_j-\delta_i)^{2}\right)^{\frac{1+\alpha}{2}}\leq\frac{dV_3(t)}{dt}\leq-\frac{1}{2}\left(\sum_{i,j=1}^na_{ij}^{\frac{2}{1+\alpha}}
(\delta_j-\delta_i)^{2}\right)^{\frac{1+\alpha}{2}}.
\end{equation}
 Moreover, since $\sum_{i,j=1}^na_{ij}^{\frac{2}{1+\alpha}}(\delta_j-\delta_i)^2=2\de^TL(\B)\de$,
 $V_3(t)=\frac{1}{2}\de^T\de$, and $\de\bot \yi$, by
Lemma \ref{yl2},
\begin{equation}\label{eqdlcomparea-2-2}
    4\lambda_2(L(\B))V_3(t)\leq\sum_{i,j=1}^na_{ij}^{\frac{2}{1+\alpha}}
(\delta_j-\delta_i)^{2}\leq 4\lambda_n(L(\B))V_3(t).
\end{equation}
By \eqref{eqdlcomparea-2-1} and \eqref{eqdlcomparea-2-2}, we have
that
\begin{equation}\label{eqdlcomparea-2-3}
    -\frac{n^{1-\alpha}}{2}\left(4\lambda_n(L(\B))V_3(t)\right)^{\frac{1+\alpha}{2}}\leq\frac{dV_3(t)}{dt}
    \leq
    -\frac{1}{2}\left(4\lambda_2(L(\B))V_3(t)\right)^{\frac{1+\alpha}{2}}.
\end{equation}

If $V_3(t)\geq\epsilon^*$, then
\begin{align*}
    \left.\frac{dV_3(t)}{dt}\right|_{\alpha=\alpha^*}\leq&-\frac{1}{2}\left(4\lambda_2(L(\B|_{\alpha=\alpha^*}))\right)^{\frac{1+\alpha^*}{2}}
    V_3(t)^{\frac{1+\alpha_*}{2}}V_3(t)^{\frac{\alpha^*-\alpha_*}{2}}
    \mbox{\ ( by \eqref{eqdlcomparea-2-3})}\\
    \leq&-\frac{1}{2}\left(4\lambda_2(L(\B|_{\alpha=\alpha^*}))\right)^{\frac{1+\alpha^*}{2}}
    n^{1-\alpha_*}2^{\alpha_*-\alpha^*}\left(\lambda_n(L(\B|_{\alpha=\alpha_*}))\right)^{\frac{1+\alpha_*}{2}}\\
    &\times\left(\lambda_2(L(\B|_{\alpha=\alpha^*}))\right)^{\frac{-(1+\alpha^*)}{2}}V_3(t)^{\frac{1+\alpha_*}{2}}\mbox{\ (by $V_3(t)\geq
    \epsilon^*$)}\\
    =&-\frac{n^{1-\alpha_*}}{2}2^{1+\alpha_*}\left(\lambda_n(L(\B|_{\alpha=\alpha_*}))\right)^{\frac{(1+\alpha_*)}{2}}
    V_3(t)^{\frac{1+\alpha_*}{2}}\\
    \leq&\left.\frac{dV_3(t)}{dt}\right|_{\alpha=\alpha_*}\mbox{ (by \eqref{eqdlcomparea-2-3}). }
\end{align*}
The case when $V_{3}(t)\leq\epsilon_*$ can be proved in the same
way.
\end{proof}
\subsection{Networks With Switching Topology}

In practice, the information channel between any two agents may not
be always available because of the restrictions of physical
equipments or the  interference in signals from external, such as
exceeding the sensing range or existence of  obstacles between
agents. Therefore, it is more reasonable to assume that the
communication topology is dynamically changing.

In order to describe the switching property of the topology, suppose
that $a_{ij}$ is a piece-wise constant right continuous function of
time, denoted by $a_{ij}(t)$, and takes value in a finite set, such
that  $a_{ii}(t)=0$ for all $i,t$. And  the changing topology is
represented  by $\mathcal{G}(\A(t))$. We only study the validity of
protocol \eqref{pro2} when the underlying topology is switching. As
for protocol \eqref{pro1}, the finite-time convergence analysis is
more challenging and we leave it as one future research topic. In
fact, for protocol \eqref{pro2}, the parameter $\alpha_{ij}$ can
also
 be changing to reflect the reliability of
information channel $(v_j,v_i)$  as weight $a_{ij}$ or to get
shorter convergence time as indicated in Theorem \ref{dlcomparea-2}.

By investigating the properties of the two protocols, we can see
that the  main difference between the two protocols is that the
 Lyapunov function $V_3(\de(t))$ used in the proof of Theorem
\ref{dl2} does not depend on the network topology. This property of
$V_3(\de(t))$ makes it a possible candidate as a common Lyapunov
function for convergence analysis of the system with switching
topology.
\begin{dl}\label{dl3}
Suppose that $\mathcal{G}(\A(t))$ is undirected and connected all
the time, all $\alpha_{ij}$ are piecewise constant right continuous
functions of time, denoted by $\alpha_{ij}(t)$, and take values in a
finite set such that $\alpha_{ij}(t)=\alpha_{ji}(t)$ and
$0<\alpha_{ij}<1$ for all $i,j,t$. Then protocol \eqref{pro2} solves
the finite-time average-agreement problem.
\end{dl}
\begin{proof}
The proof is similar to that of Theorem \ref{dl2} and we take the
same notations as in the proof of Theorem \ref{dl2}, such as
$\kappa$, $\de(t)$ and $\B$. Then $\kappa$ is also time-invariant.
Let $\alpha_0=\max_{i,j,t}\alpha_{ij}(t)$ and take the Lyapunov
function $V_3(\de(t))$. At any time $t$, we have an estimation of
$K_4(t)$, which can be chosen to depend only on $\A(t)$ and the
initial state $\x(0)$. Since all possible $\A(t)$ and
$\alpha_{ij}(t)$ are finite, $K_6=\min_tK_4(t)\lambda_2(L(\B(t)))$
exists and is larger than $0$. Then
\[\frac{d V_3(t)}{d t}\leq -2^{\alpha_0} K_6^{\frac{1+\alpha_0}{2}}V_3(t)^{\frac{1+\alpha_0}{2}}.\]

Therefore, $V_3(t)$ will reach zero in finite time
$t^*=\frac{2^{1-\alpha_0}V_3(0)^\frac{1-\alpha_0}{2}}{(1-\alpha_0)K_6^{\frac{1+\alpha_0}{2}}}$
and the switching system solves the finite-time average-agreement
problem.
\end{proof}

{\it Remark:} The conditions in Theorem \ref{dl3} can be relaxed in
several ways. For example, we can assume that the sum of time
intervals, in which $\mathcal{G}(\A(t))$ is connected, is larger
than $t^*$. In fact $\G(\A(t))$ is not necessarily undirected. This
condition can be replaced by that $\G(\A(t))$ is always
detail-balanced with a common $\w>0$ such that $\diag(\w)^T\A(t)$ is
symmetric.

In \cite{A. Jadbabaie J. Lin and A. S. Morse,W. Ren and R. W. Beard}
and \cite{L. Moreau 1}, it was shown that under certain conditions,
if the union of graph $\G(\A(t))$ across a bounded time interval is
connected or has a spanning tree, the consensus problem is solvable
asymptotically. However, in Theorem \ref{dl3}, if $\G(A(t))$ is not
connected at any time, we cannot conclude that protocol \eqref{pro2}
solves a finite-time consensus problem even though the union of
$\G(\A(t))$ across the time interval with a given length is
connected.
\begin{lz}[Counterexample]
Let
\[\A_1=\left[
         \begin{array}{ccc}
           0 & 1 & 0 \\
           1 & 0 & 0 \\
           0 & 0 & 0 \\
         \end{array}
       \right],\mbox{ and }
\A_2=\left[
         \begin{array}{ccc}
           0 & 0 & 0 \\
           0 & 0 & 1\\
           0 & 1 & 0 \\
         \end{array}
       \right].
\]
 Suppose that the system consists of three agents, communication topology
 $\G(\A(t))|_{t\in[2k, 2k+1)}=\G(\A_1)$, $\G(\A(t))|_{t\in[2k+1,
 2(k+1))}=\G(\A_2)$, $k=0,1,2,\dots$, and protocol \eqref{pro2} is applied where $\alpha_{ij}=0.5$ for all $i,j$.
 Consider the case when $\G(\A_1)$ is in effect.  Obviously $x_3(t)$ is
 time-invariant and thus we only study the dynamics of agents $1$ and
 $2$ in time-interval $[2k, 2k+1)$. Let
 $\kappa=\frac{1}{2}(x_1(2k)+x_2(2k))$, $\delta_1(t)=x_1(t)-\kappa$,
 $\delta_2(t)=x_2(t)-\kappa$ and $V(t)=\frac{1}{2}(\delta_1(t)^2+\delta_2(t)^2)$.
 \[\frac{dV(t)}{dt}=-2^{1.5}V(t)^{0.75}.\]
If $V(2k)<0.25$, then $V(t)$ will reach zero before time $2k+1$. We
can get the same result for $\G(\A_2)$.
 Given initial state $\x(0)=[0.3,0.4,0]^T$, $|\delta_1(t)|\leq\max_ix_i(t)-\min_ix_i(t)\leq\max_ix_i(0)-\min_ix_i(0)=0.4$. In the same way, $|\delta_2(t)|\leq 0.4$.
 Therefore, $V(t)\leq0.16$, and thus, the states of agents $1$ and
 $2$ are the same at time $2k+1$. In the same way, the states of agents $2$ and
 $3$ are  the same at time $2, 4,6,\dots$. Therefore,
 \[\x(2k)=\left[
            \begin{array}{ccc}
              0.5&    0.5&        0\\
    0.25&    0.25& 0.5\\
    0.25&0.25&  0.5\\
            \end{array}
          \right]^k\x(0), k=1,2,\dots.
 \]
Obviously, states of agents cannot reach consensus in  finite time.
In fact, this system asymptotically solves the average-consensus
problem.

\end{lz}

\section{Simulations}
In this section, we  present some numerical simulations to
illustrate our theoretical results.
\begin{figure}[htpb]\centering
      \includegraphics[scale=1]{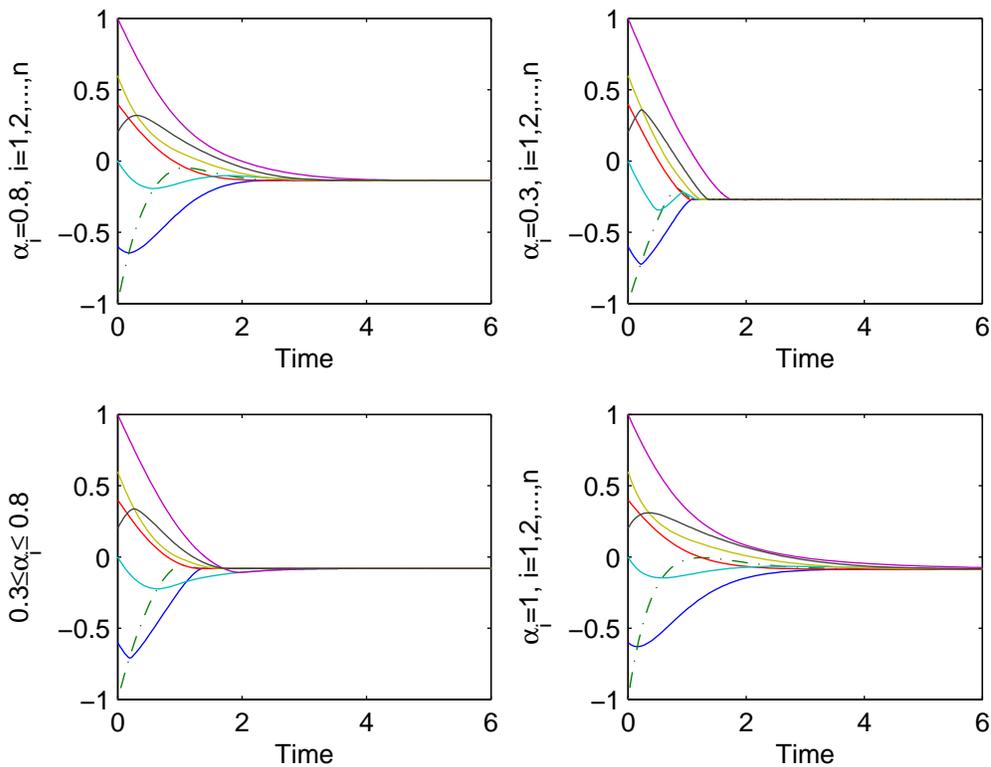}\\
      \caption{State trajectories of agents under the communication topology given in Fig. \ref{demostep} (a) when protocol
      \eqref{pro1} is applied.}
      \label{figproc1}
\end{figure}

Fig. \ref{figproc1} shows the states of agents under protocol
\eqref{pro1}. The communication topology is given in Fig.
\ref{demostep} (a), and the initial state is $[-0.6,
-1,0.4,0,1,0.6,0.2]^T$. For the case when $0.3\leq \alpha_i\leq
0.8$, we let $\alpha_1=0.3$, $\alpha_2=0.5$, $\alpha_3=0.7$,
$\alpha_4=0.8$, $\alpha_5=0.5$, $\alpha_6=0.6$ and $\alpha_7=0.55$.
As a comparison, the last one in Fig. \ref{figproc1} is the states
of agents under the typical linear consensus protocol in \cite{R.
Olfati-Saber and R. M.  Murray 1}. Since the agents' states differ a
little from each other, we can see that the system converges fast
when $\alpha_i$ is relatively small.

\begin{figure}[htpb]\centering
      \includegraphics[scale=0.6]{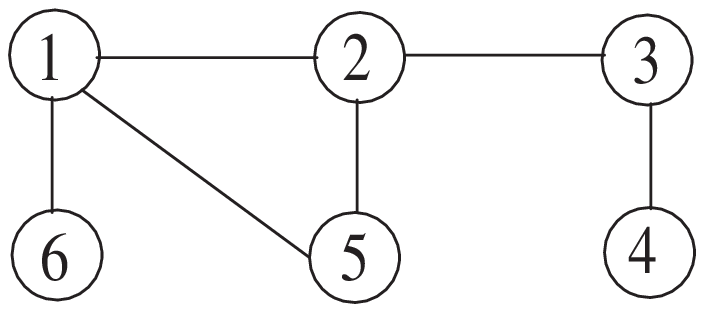}
      \includegraphics[scale=0.6]{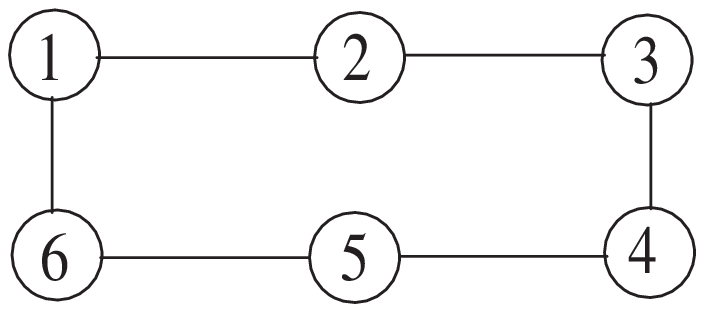}\\
      (a)\ \ \ \ \ \ \ \ \ \ \  \ \ \ \ \ \ \  \ \ \ \ \ \ \  \ \ \ (b)\\
      \includegraphics[scale=0.6]{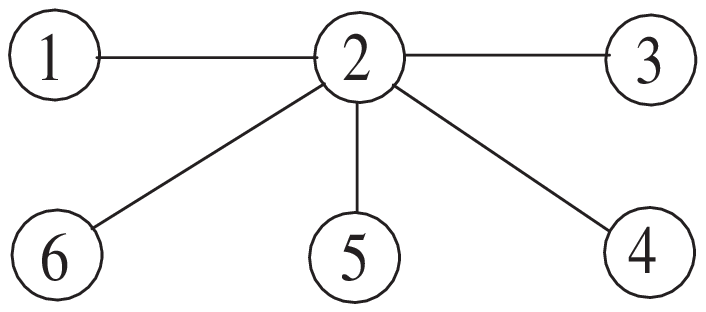}
      \includegraphics[scale=0.6]{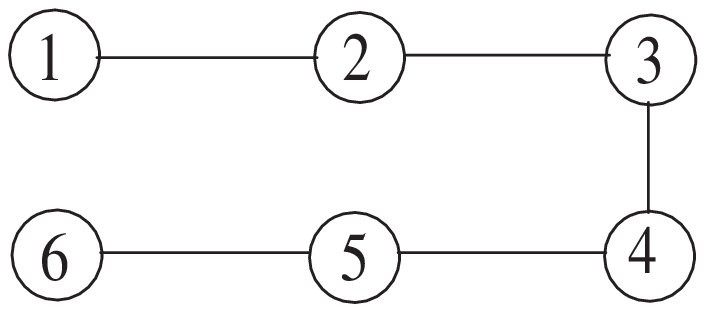}\\
        (c)\ \ \ \ \ \ \ \  \ \ \ \ \ \ \ \ \ \  \ \ \ \ \ \ \  \ \ \ (d)\\
      \caption{Four graphs: (a) $\mathcal{G}_1$, (b) $\mathcal{G}_2$, (c) $\mathcal{G}_3$, (d) $\mathcal{G}_4$}
      \label{fig1}
\end{figure}

The following simulations are performed with six agents. Fig.
\ref{fig1} shows four connected undirected  graphs and the weight of
each edge is $2$.

\begin{figure}[htpb]\centering
      \includegraphics[scale=0.8]{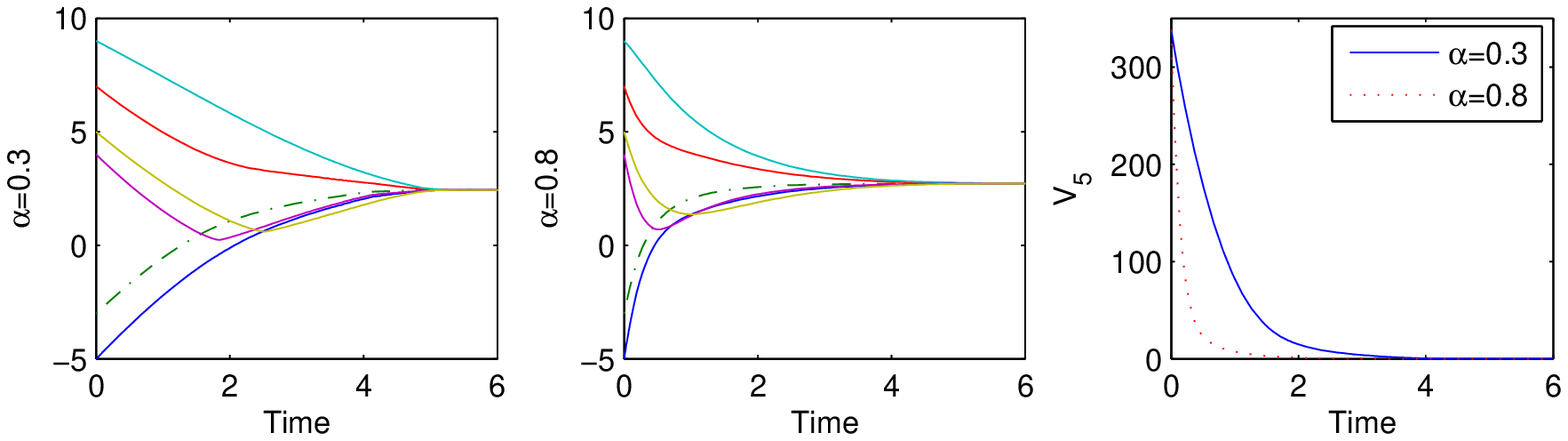}\\
      \includegraphics[scale=0.8]{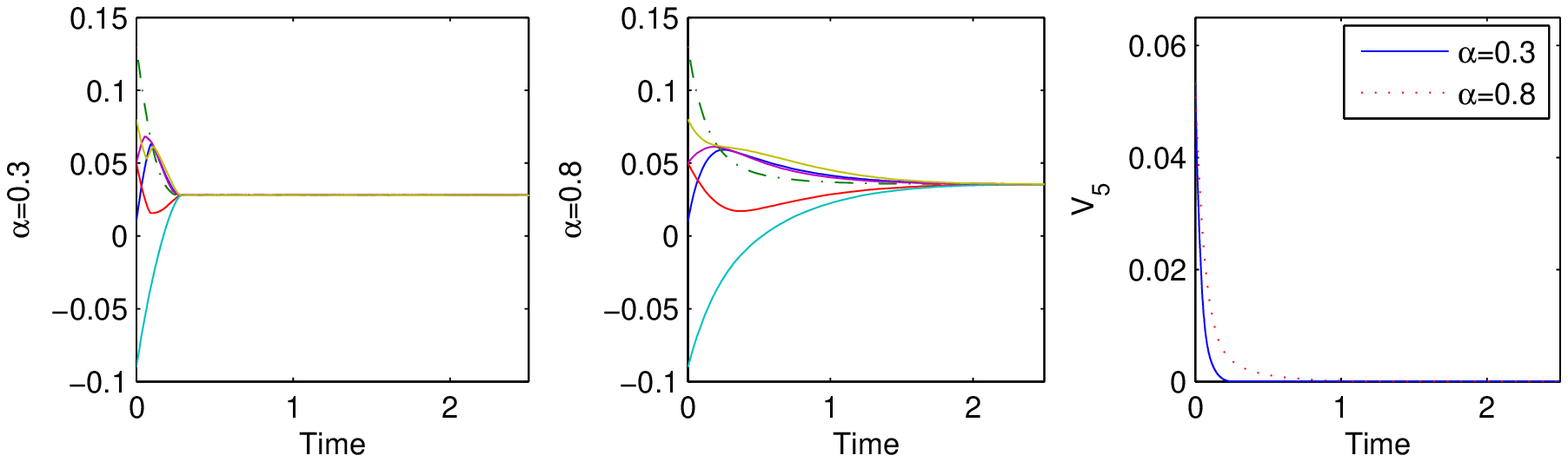}\\
      \caption{Comparison of convergence rates of the two systems under communication
      topology
       $\mathcal{G}_1$, where protocol
      \eqref{pro1} is applied, and the two systems are with different protocol parameters $\alpha_i$.}
      \label{figcmpraproc1}
\end{figure}

\begin{figure}[htpb]\centering
      \includegraphics[scale=0.8]{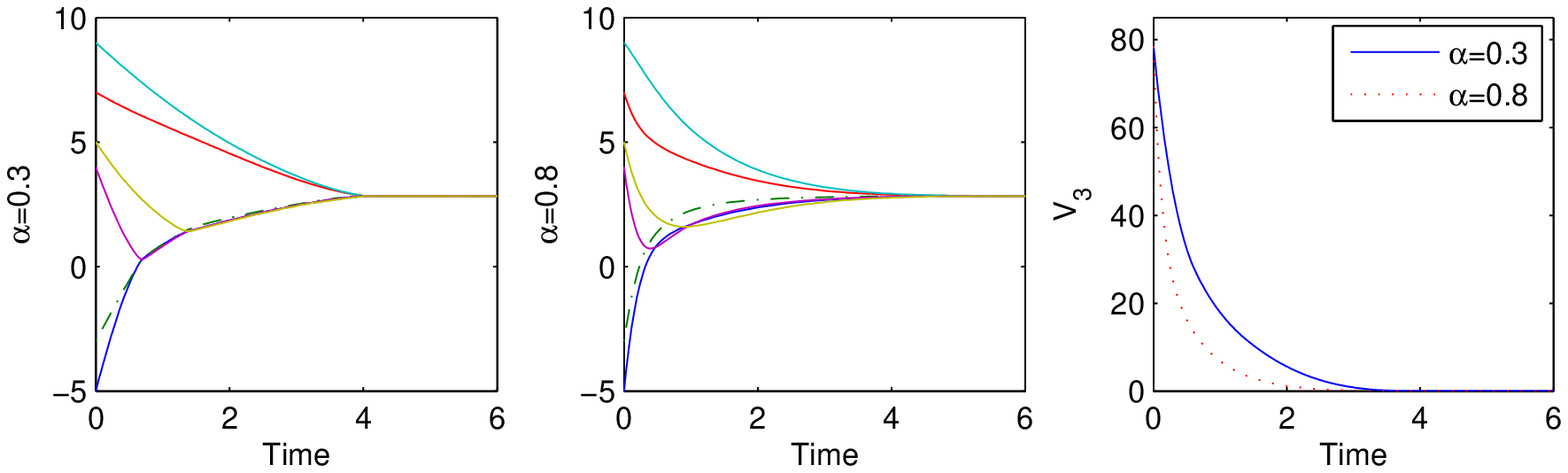}\\
      \includegraphics[scale=0.8]{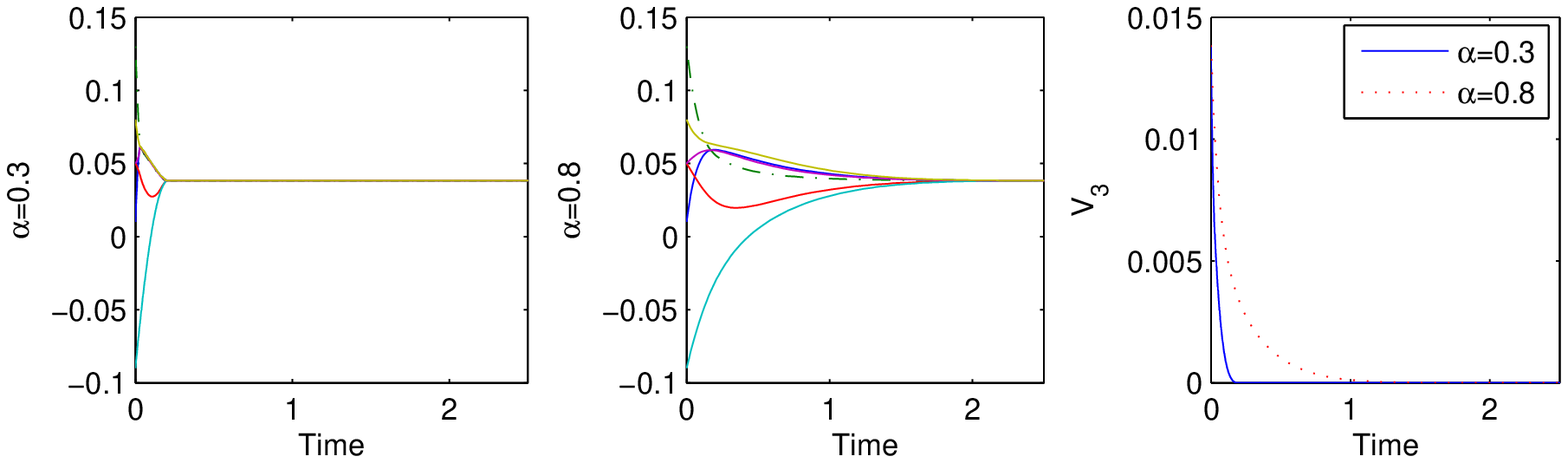}\\
      \caption{Comparison of convergence rates of the two systems under communication
      topology
       $\mathcal{G}_1$, where protocol
      \eqref{pro2} is applied, and the two systems are with different protocol parameters $\alpha_{ij}$.}
      \label{figcmpraproc2}
\end{figure}

To compare the convergence rates of of systems under protocol
\eqref{pro1} or protocol \eqref{pro2} with different $\alpha_{i}$ or
$\alpha_{ij}$, we assume the underlying topology is $\G_1$ and all
$\alpha_{i}$ or $\alpha_{ij}$ are equal to $\alpha$. For protocol
\eqref{pro1}, $\alpha$ takes $0.3$ and $0.8$ separately. By Theorem
\ref{dlcomparea-1}, if $V_5(t)\geq 7.4353$,
$\left.\frac{dV_5(t)}{dt}\right|_{\alpha=0.8}\leq\left.\frac{dV_5(t)}{dt}\right|_{\alpha=0.3}$,
and if $V_5(t)\leq 0.0569$,
$\left.\frac{dV_5(t)}{dt}\right|_{\alpha=0.3}\leq\left.\frac{dV_5(t)}{dt}\right|_{\alpha=0.8}$.
The initial state is chosen to be $[-5,-3,7,9,4,5]^T $
$(V_5(0)=338)$ and $[0.01,0.13,0.05,-0.09,0.05,0.08]^T$
$(V_5(0)=0.0533)$ respectively. The states of agents are shown in
Fig \ref{figcmpraproc1}. For protocol \eqref{pro2}, $\alpha$ also
takes $0.3$ and $0.8$ separately. By Theorem \ref{dlcomparea-2}, if
$V_3(t)\geq 42674$,
$\left.\frac{dV_3(t)}{dt}\right|_{\alpha=0.8}\leq\left.\frac{dV_3(t)}{dt}\right|_{\alpha=0.3}$,
and if $V_3(t)\leq 0.000029$,
$\left.\frac{dV_3(t)}{dt}\right|_{\alpha=0.3}\leq\left.\frac{dV_3(t)}{dt}\right|_{\alpha=0.8}$.
In fact, $\epsilon^*$ or $\epsilon_*$ given in Theorem
\ref{dlcomparea-2} is too large or small, and it can be chosen to be
 smaller or larger.  In our example, the initial state is also
set to be $[-5,-3,7,9,4,5]^T $ $(V_3(0)=78.4167)$ and
$[0.01,0.13,0.05,-0.09,0.05,0.08]^T$ $(V_3(0)=0.0138)$ respectively.
The states of agents are shown in Fig \ref{figcmpraproc2}.

\begin{figure}[htpb]\centering
      \includegraphics[scale=1]{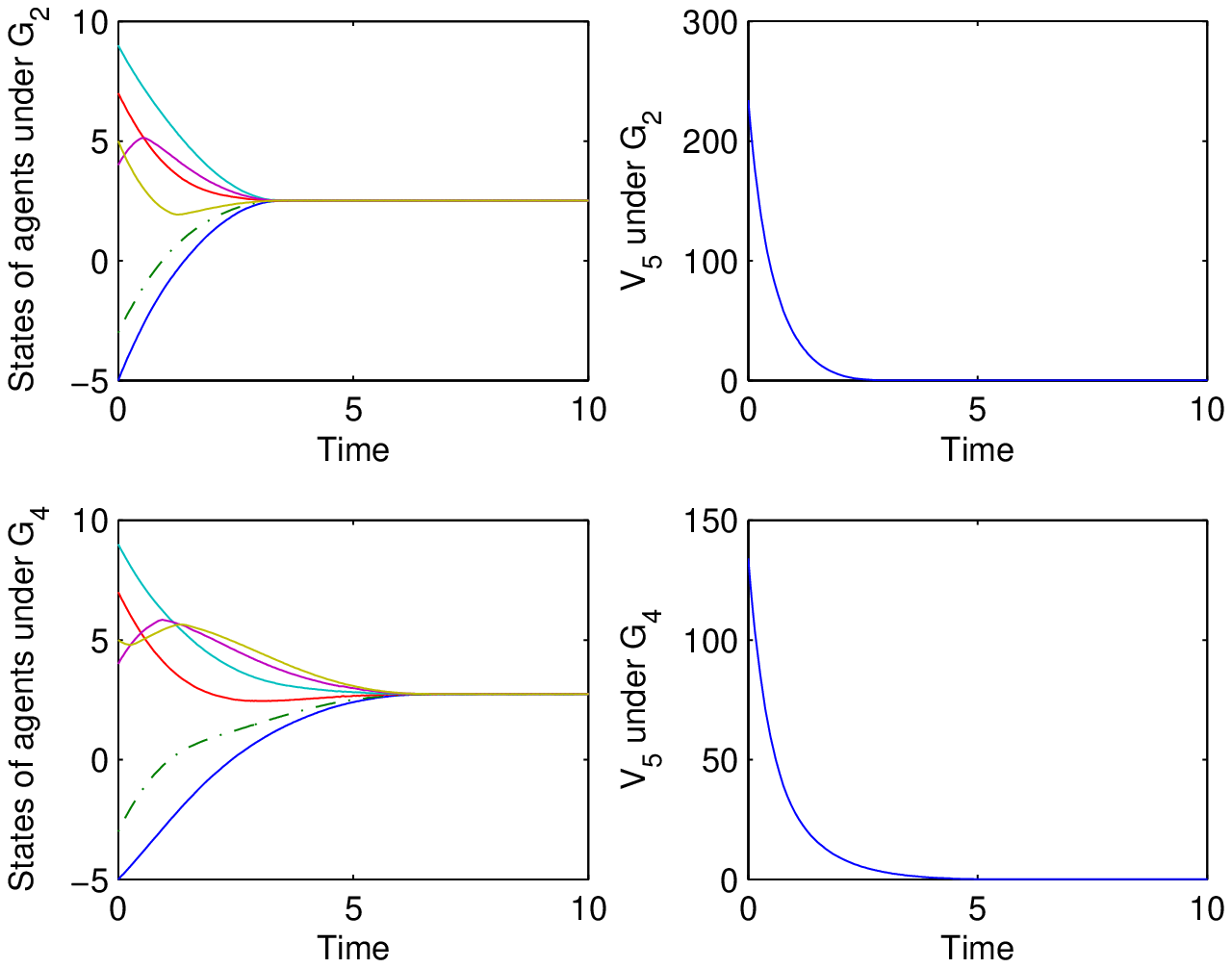}\\
      \caption{State trajectories of agents under communication topologies $\mathcal{G}_2$ and $\mathcal{G}_4$
      when protocol
      \eqref{pro1} is applied}
      \label{figproc1g2g4}
\end{figure}

\begin{figure}[htpb]\centering
      \includegraphics[scale=1]{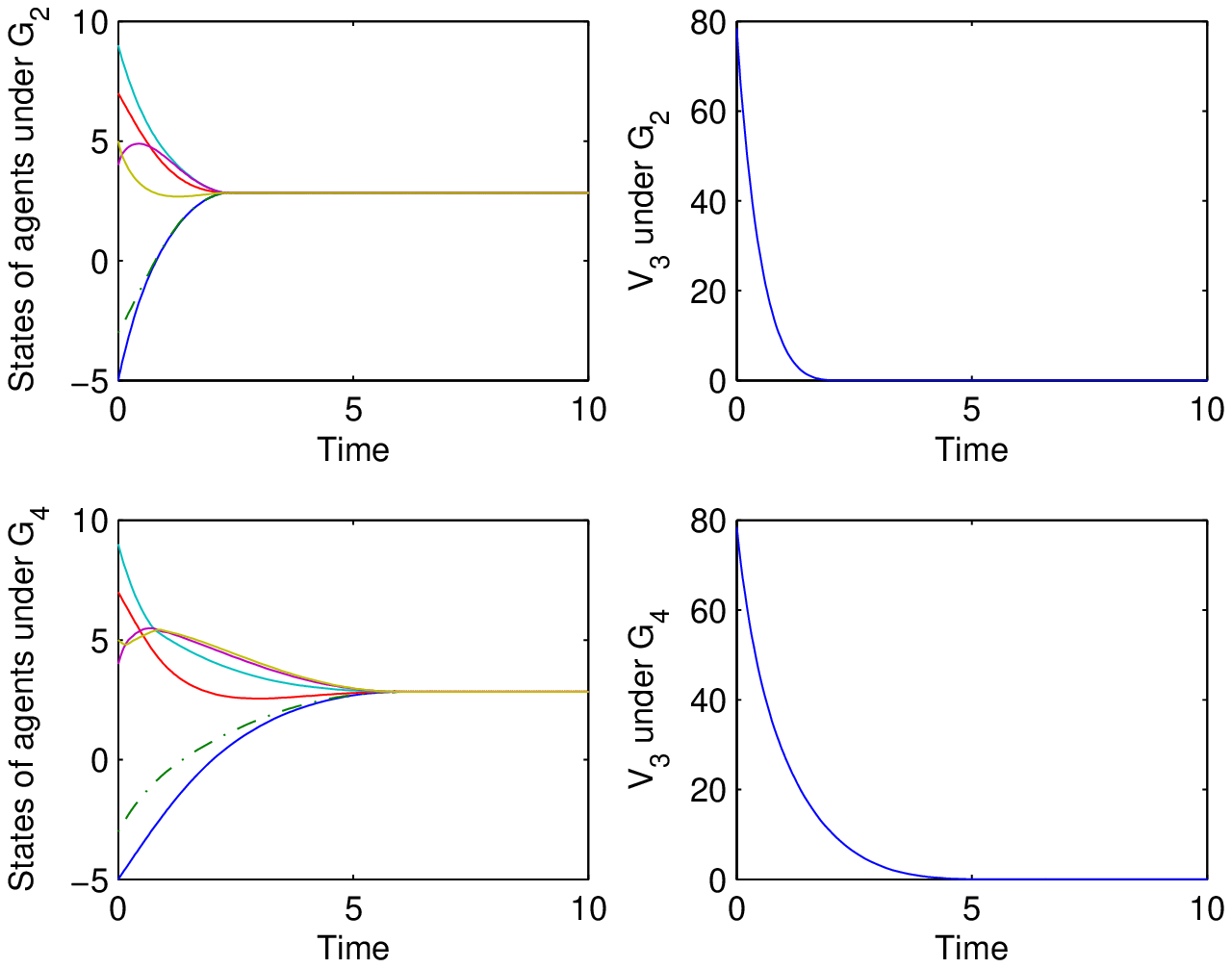}\\
      \caption{State trajectories of agents under communication topologies $\mathcal{G}_2$ and
      $\mathcal{G}_4$ when protocol
      \eqref{pro2} is applied}
      \label{figproc2g2g4}
\end{figure}

To show the effect of algebraic connectivity on the convergence
time, consider the case when communication topology is $\G_2$ and
$\G_4$ separately and $\alpha_{ij}$ and $\alpha_i$ are all equal to
$\alpha=0.5$. Let $\x(0)=[-5,-3,7,9,4,5]^T$. The algebraic
connectivity of $\G_2$ is $2$ and that of $\G_4$ is $0.5359$. The
upper bounds of convergence times estimated by differential equation
\eqref{ineqv5-a} are $6.0638$ and $16.2819$ separately. And the
algebraic connectivity corresponding to
$\G([a_{ij}^{\frac{2}{1+\alpha}}])$ of $\G_2$ is $2.5198$, and that
of $\G_4$ is $0.6752$. The estimated upper bounds of convergence
times by differential equation \eqref{ineqv3-a} are $4.2085$ and
$11.2999$. The states of agents under protocols \eqref{pro1} and
\eqref{pro2} are shown in Fig. \ref{figproc1g2g4} and Fig.
\ref{figproc2g2g4} separately.

\begin{figure}[htpb]\centering
      \includegraphics[scale=1]{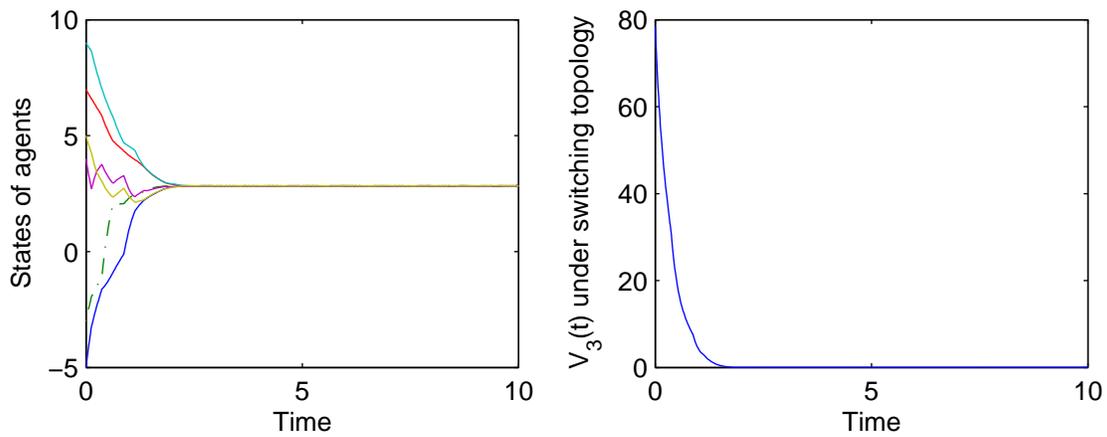}\\
      \caption{State trajectories of agents under switching communication topology}
      \label{figswitching}
\end{figure}

If the communication topology is switching from $\mathcal{G}_1$, to
$\mathcal{G}_2$, to $\mathcal{G}_3$, to $\mathcal{G}_4$, and back to
$\mathcal{G}_1$, periodically, and each of them lasts for $0.25$
seconds, then we apply protocol \eqref{pro2} and the states of
agents achieve the average-agreement in finite time, the estimated
upper bound of which is $t_3=11.2999$. The trajectories of them are
shown in Fig. \ref{figswitching}.

\section{Conclusion}
We have discussed the finite-time consensus problems for multi-agent
systems and presented two effective continuous finite-time consensus
protocols. Furthermore, the relationship between the convergence
time and the communication topology, the initial states and the
protocol parameters $\alpha_i$ or $\alpha_{ij}$ was analyzed.
Several simulations demonstrated the effectiveness of our
theoretical results.

The work of this paper is the first step toward the finite-time
consensus analysis of multi-agent systems, and there are still some
other interesting and important topics to be addressed. For example,
does  protocol \eqref{pro2} still work when the communication
topology has a spanning tree? If the system is under switching
topology and with communication time-delays, do there exist similar
results? Do there exist other effective finite-time protocols? These
problems are currently under investigations.

\appendices
\section{}
{\it Ger\v sgorin Disk Theorem (\cite{Horn}, pp. 344, Theorem
6.1.1):} Let $\A\in[a_{ij}]\in\C^{n\times n}$, and let
\[R'_i(\A)=\sum_{j=1 \atop i\not=j}^n|a_{ij}|, 1\leq i\leq n\]
denote the deleted absolute row sums of $\A$. Then all the
eigenvalues of $\A$ are located in the union of $n$ discs
\[
\bigcup_{i=1}^n\{c\in\C:|c-a_{ii}|\leq R'_i(\A)\}.
\]

\begin{yl}[see \cite{Horn}, Theorem 6.2.24]\label{ylHorn6224}
Let  $\A\in \R^{n\times n}$ be a  nonnegative matrix. The following
are equivalent:
\begin{enumerate}
    \item $\A$ is irreducible;
    \item $(\bI+\A)^{n-1}>0$;
    \item $\G(\A)$ is strongly connected.
\end{enumerate}
\end{yl}

{\it Perron-Frobenius Theorem (\cite{Horn}, pp. 508, Theorem
8.4.4):} Let $\A\in\R^{n\times n}$ and suppose that $\A$ is
irreducible and nonnegative. Then
\begin{enumerate}
    \item $\rho(\A)>0$;
    \item $\rho(\A)$ is an eigenvalue of $\A$;
    \item There is a positive vector $\x$ such that
    $\A\x=\rho(\A)\x$; and
    \item $\rho(\A)$ is an algebraically (and hence geometrically)
    simple eigenvalue of $\A$.
\end{enumerate}

{\it Peano's Existence Theorem (\cite{P. Hartman}, pp. 10, Theorem
2.1):} Let $\y, f\in\R^n$; $f(t,\y)$ continuous on $R:t_0\leq t\leq
t_0+a$, $\|\y-\y_0\|_{\infty}\leq b$; $M$ a bound for
$\|f(t,\y)\|_{\infty}$ on $R$; $\alpha=\min(a,b/M)$. Then
differential equation
\[\frac{d\y}{dt}=f(t,\y), \y(t_0)=\y_0\]
possesses at least one solution $\y=\y(t)$ on $[t_0, t_0+\alpha]$.

{\it Extension Theorem (\cite{P. Hartman}, pp. 12, Theorem 3.1):}
 Let $f(t,\y)$ be continuous on an open
$(t,\y)$-set $E$ and let $\y(t)$ be a solution of $\y'=f(t,\y)$ on
some interval. Then $\y(t)$ can be extended (as a solution) over a
maximal interval of existence $(\omega_-,\omega_+)$. Also, if
$(\omega_-,\omega_+)$ is a maximal interval of existence, then
$\y(t)$ tends to the boundary $\partial E$ of $E$ as $t\to\omega_-$
and $t\to\omega_+$.

{\it Comparison Principle (\cite{P. Hartman}, pp. 26, Theorem 4.1):}
Let $U(t,u)$ be continuous on an open $(t,u)$-set $E$ and $u=u^0(t)$
the maximal solution of
\[\frac{du}{dt}=U(t,u), u(t_0)=u_0.\]
Let $v(t)$ be a continuous function on $[t_0, t_0 +a]$ satisfying
the conditions $u(t_0)\leq u_0$, $(t,v(t))\in E$, and $v(t)$ has a
right derivative $D_R(v(t))$ on $t_0\leq t<t_0+a$ such that
\[D_Rv(t)\leq U(t,v(t)).\]
Then, on a common interval of existence of $u^0(t)$ and $v(t)$,
\[v(t)\leq u^0(t).\]

\end{document}